\documentclass{article} 
\usepackage{latexsym, a4wide}
\usepackage{amsmath, rotating, color} 
\usepackage{amsfonts,amssymb, amsthm, mathtools} 
\usepackage{pictexwd} 
\usepackage{mathrsfs}
\usepackage[normalem]{ulem} 
\newcommand\unnumberedfootnote[1]{ %
  \let\temp=\thefootnote %
  \renewcommand{\thefootnote}{}%
  \footnote{#1}%
  \let\thefootnote=\temp%
  \addtocounter{footnote}{-1}}

\newtheorem{theorem}{Theorem}
\newtheorem{proposition}{Proposition}[section]
\newtheorem{lemma}[proposition]{Lemma}

\theoremstyle{definition}
\newtheorem{remark}[proposition]{Remark}

\newtheorem{step}{Step}

\numberwithin{equation}{section}
\begin{document}
\title{\LARGE The tree length of an evolving coalescent}

\thispagestyle{empty}

\author{{\sc by P. Pfaffelhuber\thanks{Corresponding Author;
      Albert-Ludwigs Universit\"at; Abteilung Mathematische
      Stochastik, Eckerstr. 1, D--79104 Freiburg; email:
      p.p@stochastik.uni-freiburg.de; Tel.: +49 761 203 5667; Fax: +49
      761 203 5661}, A. Wakolbinger\thanks{Goethe-Universit\"at;
      Institut f\"ur Mathematik; Robert-Mayer-Str. 10; D--60325
      Frankfurt; email: wakolbinger@math.uni-frankfurt.de; Tel: +49 69
      798 28651; Fax: +49 69 798 28841} and
    H. Weisshaupt\thanks{Albert-Ludwigs Universit\"at; Center for
      Biosystems Analysis; Habsburger Str. 49; email:
      heinz.weisshaupt@zbsa.de;
      Tel: +49 761 203 97170; Tel: +49 761 203 97217}}\\[2ex]
  \date{}}

\maketitle
\unnumberedfootnote{\emph{AMS 2000 subject classification.} 60K35
(Primary) 92D25 (Secondary).}

\unnumberedfootnote{\emph{Keywords and phrases.}  Kingman's
  Coalescent, genealogical trees, Moran model, evolution of tree
  length, large population limit, Gumbel distribution}

\begin{abstract}
\noindent
A well-established model for the genealogy of a large population in
equilibrium is Kingman's coalescent.  For the population together with
its genealogy evolving in time, this gives rise to a time-stationary
tree-valued process. We study the sum of the branch lengths, briefly
denoted as tree length, and prove that the (suitably compensated)
sequence of tree length processes converges, as the population size
tends to infinity, to a limit process with c{\`a}dl\`ag paths,
infinite infinitesimal variance, and a Gumbel distribution as its
equilibrium.
\end{abstract}

\section{Introduction}
Kingman's coalescent \cite{Kingman1982a, Kingman1982} is a widely used
model for the single-locus genealogy in a population, see
\cite{Wakeley2008} and references therein. It arises in a suitable
rescaling of time under the assumptions of a neutral evolution and an
exchangeable reproduction dynamics with short-tailed offspring
distribution. An intuitive way to think of Kingman's coalescent is to
imagine a random tree with infinitely many leaves at time $t$, where
backwards in time any two lineages independently coalesce at rate
1. Taking $\mathrm N$ instead of infinitely many leaves gives
Kingman's $\mathrm N$-coalescent. The latter figures as the genealogy
of an $\mathrm N$-sample taken from a large population, and also as
the genealogy of the total population in a standard Moran model with
population size $\mathrm N$.

Two functionals of coalescent trees are of particular interest: the
distance from the root to the leaves, or \emph{depth}, and the sum of
branch lengths, or \emph{tree length}. It is well known that the
expected depth of Kingman's $\mathrm N$-coalescent equals
$2(1-\frac{1}{\mathrm N})$, whereas its expected tree length is $\sim
2\log \mathrm N$ as $\mathrm N\to \infty$.  More can be said: when
compensated by $2\log \mathrm N$, half the tree length of Kingman's
$\mathrm N$-coalescent converges in law to a Gumbel distributed random
variable (having the cumulative distribution function $x\mapsto
e^{-e^{-x}}$). This result can be read off from \cite[p. 153, first
equation]{Tavare1984}; see also \cite{WiufHein1999},
\cite{Tavare2004}, \cite{DrmotaEtAl2007} and \cite{Wakeley2008}.

With individual offspring distributions that are not short-tailed,
coalescents different from Kingman's appear as the genealogies of
large populations.  In the so-called $\Lambda$-coalescents
\cite{Pitman1999}, more than two lines can coalesce, giving rise to
{\em multiple mergers}, and asymptotic tree length distributions arise
that are different from Gumbel distributions. For special classes
including those of Beta-coalescents, results on the asymptotic tree
length were obtained in \cite{Moehle2006}, \cite{DrmotaEtAl2007},
\cite{BBS2009} and \cite{DelmasEtAl2008}.


With a population evolving in time, its genealogical relationships
evolve as well. Their evolution is described by a \emph{tree-valued
  process} \cite{GrevenPfaffelhuberWinter2010} similarly as the
change of allele frequencies is captured by measure-valued diffusions
\cite{Dawson1993, Etheridge2001}. Jumps of the tree depth correspond
to the loss of one of the currently two oldest families from the
population, and hence to the establishment of a new most recent common
ancestor (MRCA) of the population. The resulting tree depth process in
the case of Kingman's coalescent was analyzed in
\cite{PfaffelhuberWakolbinger2006} and \cite{DelmasEtAl2010}.

In the present paper we focus on the (compensated) tree length in
Kingman's coalescent and describe its evolution in an infinite
population. Our main result is that this process has c\`adl\`ag paths
and infinite infinitesimal variance (Theorem~\ref{TH1}).  As already
stated, the one-dimensional projections of this process are Gumbel
distributed.

We construct the process of compensated tree length, denoted by
$\mathscr L$, as a limit using tools from weak convergence of
processes. In addition, we also provide a strong convergence result,
i.e.\ a version of Theorem \ref{TH1} in terms of convergence towards
$\mathscr L$ in probability. For this, we use the lookdown process
introduced in \cite{DonnellyKurtz1999}, which provides genealogies of
Moran models of any population size on one and the same probability
space. Our Proposition \ref{P:TH1} shows that on this space, the
compensated Kingman tree lengths lead to a c\`adl\`ag path-valued
limit in probability.
Hence, the process $\mathscr L$ can be defined directly in terms of a
sequence Moran models -- or in terms of the lookdown graph -- and as
such is a natural object to study. Some challenging questions remain,
e.g.  a) Is the limit robust in the sense that $\mathscr L$ describes
also the limiting tree length process for (a large class of) Cannings
models with short-tailed offspring distributions?  b) Is there an
intrinsic characterization of $\mathscr L$ in terms of a stochastic
dynamics? In particular, is $\mathscr L$ a semimartingale?

~

The length of a coalescent is of relevance in empirical population
genetics \cite{Wakeley2008}. In the infinite sites model the number
of mutations seen in a population of size $\mathrm N$ at time $t$ is
Poisson distributed with parameter proportional to tree length and to
the mutation rate.  The process of tree lengths has also attracted
interest in the study of diversity in real populations
\cite[Fig.~2c]{RauchBarYam2004}. There, sudden losses of diversity in
a population are related to jumps of the tree length process. Such
jumps occur at any resampling event and correspond to the length of an
external branch breaking off the tree. The asymptotics of external
branch lengths are investigated in detail in
\cite{CaliebeNeiningerKrawczakRoesler2007}, see also Remark \ref{ebl}
and Section \ref{S:ext}.

~

Our paper is organized as follows. After specifying the model we
present our results on weak (Section \ref{S:res}) and strong
(Section~\ref{S:strong}) convergence of tree lengths and tree length
processes. In Section~\ref{S:3} we provide some auxiliary results on
Kingman's coalescent for fixed times and on Moran models. Section
\ref{S:Sec4} completes the proof of Theorem \ref{TH1}, and Section
\ref{S:proofs2} contains the proofs of the strong convergence results,
Propositions \ref{L:Gumbel2} and \ref{P:TH1}.

\section{Convergence of tree length distributions}
\label{S:res}
\begin{figure}
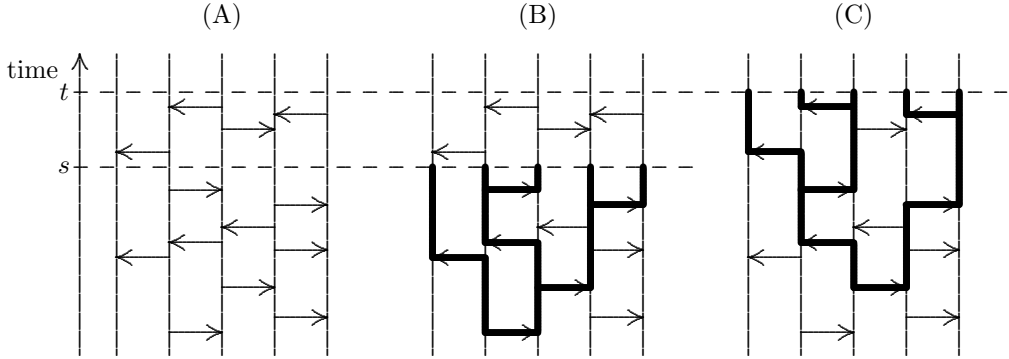

\beginpicture
\setcoordinatesystem units <.7cm,1cm>
\setplotarea x from 0.5 to 18, y from 1 to 5.5
\put{(A)} [cC] at 3 5.5
\put{(B)} [cC] at 9 5.5
\put{(C)} [cC] at 15 5.5

\plot 1 1 1 5 /
\plot 2 1 2 5 /
\plot 3 1 3 5 /
\plot 4 1 4 5 /
\plot 5 1 5 5 /
\arrow <0.2cm> [0.375,1] from 2 1.3 to 3 1.3
\arrow <0.2cm> [0.375,1] from 4 1.5 to 5 1.5
\arrow <0.2cm> [0.375,1] from 3 1.9 to 4 1.9
\arrow <0.2cm> [0.375,1] from 2 2.3 to 1 2.3
\arrow <0.2cm> [0.375,1] from 4 2.4 to 5 2.4
\arrow <0.2cm> [0.375,1] from 3 2.5 to 2 2.5
\arrow <0.2cm> [0.375,1] from 4 2.7 to 3 2.7
\arrow <0.2cm> [0.375,1] from 4 3.0 to 5 3.0
\arrow <0.2cm> [0.375,1] from 2 3.2 to 3 3.2
\arrow <0.2cm> [0.375,1] from 2 3.7 to 1 3.7
\arrow <0.2cm> [0.375,1] from 3 4.0 to 4 4.0
\arrow <0.2cm> [0.375,1] from 5 4.2 to 4 4.2
\arrow <0.2cm> [0.375,1] from 3 4.3 to 2 4.3

\plot 7 1 7 5 /
\plot 8 1 8 5 /
\plot 9 1 9 5 /
\plot 10 1 10 5 /
\plot 11 1 11 5 /
\arrow <0.2cm> [0.375,1] from 8 1.3 to 9 1.3
\arrow <0.2cm> [0.375,1] from 10 1.5 to 11 1.5
\arrow <0.2cm> [0.375,1] from 9 1.9 to 10 1.9
\arrow <0.2cm> [0.375,1] from 8 2.3 to 7 2.3
\arrow <0.2cm> [0.375,1] from 10 2.4 to 11 2.4
\arrow <0.2cm> [0.375,1] from 9 2.5 to 8 2.5
\arrow <0.2cm> [0.375,1] from 10 2.7 to 9 2.7
\arrow <0.2cm> [0.375,1] from 10 3.0 to 11 3.0
\arrow <0.2cm> [0.375,1] from 8 3.2 to 9 3.2
\arrow <0.2cm> [0.375,1] from 8 3.7 to 7 3.7
\arrow <0.2cm> [0.375,1] from 9 4.0 to 10 4.0
\arrow <0.2cm> [0.375,1] from 11 4.2 to 10 4.2
\arrow <0.2cm> [0.375,1] from 9 4.3 to 8 4.3

\multiput {{\color{black}{\tiny $\bullet$}}} at 7 3.5 *120 0 -0.01 /
\multiput {{\color{black}{\tiny $\bullet$}}} at 7 2.3 *100 0.01 0 /
\multiput {{\color{black}{\tiny $\bullet$}}} at 8 2.3 *100 0 -0.01 /
\multiput {{\color{black}{\tiny $\bullet$}}} at 8 1.3 *100 0.01 0 /
\multiput {{\color{black}{\tiny $\bullet$}}} at 8 3.5 *100 0 -0.01 /
\multiput {{\color{black}{\tiny $\bullet$}}} at 8 2.5 *100 0.01 0 /
\multiput {{\color{black}{\tiny $\bullet$}}} at 9 2.5 *120 0 -0.01 /
\multiput {{\color{black}{\tiny $\bullet$}}} at 9 3.5 *30 0 -0.01 /
\multiput {{\color{black}{\tiny $\bullet$}}} at 9 3.2 *100 -0.01 0 /
\multiput {{\color{black}{\tiny $\bullet$}}} at 10 3.5 *160 0 -0.01 /
\multiput {{\color{black}{\tiny $\bullet$}}} at 10 1.9 *100 -0.01 0 /
\multiput {{\color{black}{\tiny $\bullet$}}} at 11 3.5 *50 0 -0.01 /
\multiput {{\color{black}{\tiny $\bullet$}}} at 11 3.0 *100 -0.01 0 /

\plot 13 1 13 5 /
\plot 14 1 14 5 /
\plot 15 1 15 5 /
\plot 16 1 16 5 /
\plot 17 1 17 5 /

\arrow <0.2cm> [0.375,1] from 14 1.3 to 15 1.3
\arrow <0.2cm> [0.375,1] from 16 1.5 to 17 1.5
\arrow <0.2cm> [0.375,1] from 15 1.9 to 16 1.9
\arrow <0.2cm> [0.375,1] from 14 2.3 to 13 2.3
\arrow <0.2cm> [0.375,1] from 16 2.4 to 17 2.4
\arrow <0.2cm> [0.375,1] from 15 2.5 to 14 2.5
\arrow <0.2cm> [0.375,1] from 16 2.7 to 15 2.7
\arrow <0.2cm> [0.375,1] from 16 3.0 to 17 3.0
\arrow <0.2cm> [0.375,1] from 14 3.2 to 15 3.2
\arrow <0.2cm> [0.375,1] from 14 3.7 to 13 3.7
\arrow <0.2cm> [0.375,1] from 15 4.0 to 16 4.0
\arrow <0.2cm> [0.375,1] from 17 4.2 to 16 4.2
\arrow <0.2cm> [0.375,1] from 15 4.3 to 14 4.3

\multiput {{\color{black}{\tiny $\bullet$}}} at 13 4.5 *80 0 -0.01 /
\multiput {{\color{black}{\tiny $\bullet$}}} at 13 3.7 *100 0.01 0 /
\multiput {{\color{black}{\tiny $\bullet$}}} at 14 3.7 *120 0 -0.01 /
\multiput {{\color{black}{\tiny $\bullet$}}} at 14 2.5 *100 0.01 0 /
\multiput {{\color{black}{\tiny $\bullet$}}} at 15 2.5 *60 0 -0.01 /
\multiput {{\color{black}{\tiny $\bullet$}}} at 15 1.9 *100 0.01 0 /
\multiput {{\color{black}{\tiny $\bullet$}}} at 14 4.5 *20 0 -0.01 /
\multiput {{\color{black}{\tiny $\bullet$}}} at 14 4.3 *100 0.01 0 /
\multiput {{\color{black}{\tiny $\bullet$}}} at 15 4.3 *110 0 -0.01 /
\multiput {{\color{black}{\tiny $\bullet$}}} at 15 3.2 *100 -0.01 0 /
\multiput {{\color{black}{\tiny $\bullet$}}} at 15 4.5 *20 0 -0.01 /
\multiput {{\color{black}{\tiny $\bullet$}}} at 16 4.5 *30 0 -0.01 /
\multiput {{\color{black}{\tiny $\bullet$}}} at 16 4.2 *100 0.01 0 /
\multiput {{\color{black}{\tiny $\bullet$}}} at 17 4.2 *120 0 -0.01 /
\multiput {{\color{black}{\tiny $\bullet$}}} at 17 3.0 *100 -0.01 0 /
\multiput {{\color{black}{\tiny $\bullet$}}} at 16 3.0 *110 0 -0.01 /
\multiput {{\color{black}{\tiny $\bullet$}}} at 17 4.5 *30 0 -0.01 /





\arrow <0.2cm> [0.375,1] from 0.3 1 to 0.3 5
\put{time} [cC] at -.6 4.8
\setdashes
\plot .2 3.5 12 3.5 /
\plot .2 4.5 18 4.5 /
\put{$s$} [cC] at 0 3.5
\put{$t$} [cC] at 0 4.5
\endpicture
\caption{\label{fig1} (A) The graphical representation of a Moran
  model of size $\mathrm{N} =5$.  Arrows between lines indicate
  resampling events.  By resampling the genealogical relationships
  between individuals change. The individual at the tip dies and the
  other one reproduces. (B) The genealogical relationships at time $s$
  can be read off from the graphical representation.  The vertical
  bold lines constitute the genealogical tree $\mathscr T^{\mathrm
    N}_s$. (C) By time $t$ the genealogical tree has changed to
  $\mathscr T^{\mathrm N}_t$. However, some parts of the tree
  $\mathscr T^{\mathrm N}_t$ were already present in the tree
  $\mathscr T^{\mathrm N}_s$.}
\end{figure}
Consider a {\em Moran model} with constant population size $\mathrm{N}
$, started at time $-\infty$.\footnote{ We use roman upper case
  letters $\mathrm K, \mathrm N,...$ for real (non-random) numbers in
  order to distinguish them from random variables denoted by $K$, $N$,
  \ldots.} Each (unordered) pair of individuals {\em resamples} at
rate 1; in any such resampling event, one of the two individuals
reproduces and the other one dies. See Figure \ref{fig1}(A) for an
illustration.

At any time $t\in\mathbb R$, the common ancestry of all individuals in
the population is described by a random genealogical tree, which is
Kingman's $\mathrm{N}$-coalescent \cite{Kingman1982a}.  With time $t$
varying, we obtain a tree-valued process denoted by $\mathscr
T^\mathrm{N} = (\mathscr T^\mathrm{N} _t)_{t\in\mathbb R}$, whose
random path we can read off from the graphical representation, see
Figure \ref{fig1}(B) and (C).

Let $\widetilde \ell$ be the map that sends a (finite) tree to its
length, i.e.\ to the sum of the lengths of all branches. Back from a
fixed time $t$, each (unordered) pair of ancestral lines coalesces at
rate 1, therefore the length of the time interval during which the
genealogical tree $\mathscr T_t^\mathrm{N}$ has $k$ lines is
exponentially distributed with the number of pairs, $\binom{k}{2}$, as
parameter. Consequently, the expectation and the variance of the tree
length are
$$ \mathbb E[\widetilde \ell(\mathscr T_t^\mathrm{N} )] = \sum_{i=2}^{\mathrm N}
i \frac{1}{\binom{i}{2}} = 2 \sum_{i=1}^{\mathrm N-1} \frac 1i, \qquad
\mathbb V[\widetilde \ell(\mathscr T_t^{\mathrm N})] =
\sum_{i=2}^{\mathrm N} i^2 \frac{1}{\binom{i}{2}^2} = 4
\sum_{i=1}^{\mathrm N-1} \frac 1{i^2}.$$ We are going to study the
{\em compensated tree length process}
\begin{align}
 \mathscr L^\mathrm{N} := (\widetilde \ell(\mathscr T^\mathrm{N}
 _t)-2\log\mathrm N)_{t\in\mathbb R}
\end{align}
in the limit $\mathrm N\to\infty$.

\begin{figure}
\mbox{}
\begin{center}
\includegraphics[width=12cm]{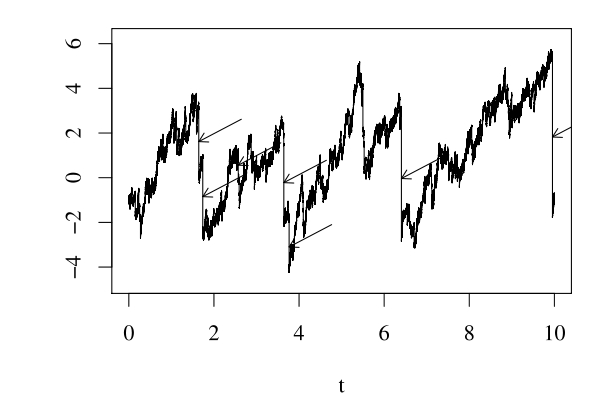}
\end{center}
\vspace{-5cm}
\hspace{2cm}
\begin{sideways} {\large$\mathscr L_t^{\mathrm 100}$} \end{sideways}
\vspace{4cm}
\caption{\label{fig:2}One realization of the tree length process
  $\mathscr L^{\mathrm 100}$ started at time $t=0$ in equilibrium. The
  arrows point to those jumps of the tree length at which the MRCA of
  the total population changes.}
\end{figure}

One realization of the process $\mathscr L^{\mathrm 100}$ can be seen
in Figure \ref{fig:2}. Here, several large jumps of the tree length
can be observed. Particularly large jumps in the tree length arise
when the MRCA of the total population changes.

\subsection{Tree lengths at fixed times}
\noindent
We recall a basic fact about the asymptotics of the law of $\mathscr
L^{\mathrm N}_t$ as $\mathrm N \to \infty$.

\begin{proposition}[Tree lengths for fixed times]\label{L:Gumbel}
  For $t\in\mathbb R$, the law of $\tfrac 12\mathscr L^{\mathrm N}_t$
  converges as $\mathrm N\to\infty$ weakly to the standard Gumbel
  distribution with cumulative distribution function $x\mapsto
  e^{-e^{-x}}$.
\end{proposition}


\begin{proof}
  We briefly repeat the argument from \cite[p. 255]{WiufHein1999}. Let
  $X_2, X_3,...$ be independent random variables such that $X_j$ has
  an exponential distribution with rate $\binom{j}{2}$. In addition,
  let $Y_1,Y_2,...$ be independent such that $Y_j$ has an exponential
  distribution with rate $j$ and $Z_1,Z_2,...$ be independent
  exponential, each with parameter $1$. Then
 $$ \tfrac 12 \widetilde\ell(\mathscr T_t^{\mathrm N})  \stackrel{d}{=} \tfrac 12
 \sum_{j=2}^{\mathrm N} j X_j \stackrel{d}{=} \sum_{j=1}^{\mathrm
   N-1} Y_j \stackrel{d}{=} \max_{1\leq j\leq \mathrm N-1} Z_j,$$
 which when shifted by $\log\mathrm N$ has the asserted limit in
 distribution as $\mathrm N\to\infty$.
\end{proof}

~

\subsection{The evolution of tree lengths}
\noindent 
Next, we come to our main result on the limit of the compensated tree
length processes $\mathscr L^{\mathrm N}$. We denote by $\mathbb D$
the space of real-valued c\`adl\`ag functions on the time axis
$(-\infty,\infty)$, equipped with the Skorokhod topology.

\begin{theorem} \label{TH1} There is a process $\mathscr L = (\mathscr
 L_t)_{t\in\mathbb R}$ with sample paths in $\mathbb D$ such that
$$ \mathscr L^\mathrm{N}  \Longrightarrow \mathscr L \text{ as }\mathrm{N} \to\infty.$$
The distribution of $\tfrac 12 \mathscr L_t$ is Gumbel for all
$t\in\mathbb R$. The process $\mathscr L$ has infinite infinitesimal
variance, with
\begin{align}\label{Thm1}
\frac{1}{t|\log t|} \mathbb E[(\mathscr L_t
- \mathscr L_0)^2] \xrightarrow{t\to 0} 4.
\end{align}
\end{theorem}

\begin{remark}[Connection to external branch lengths]\label{ebl}
  An important characteristics of the tree length process is its jump
  size statistics. For this, let $F$ be a randomly chosen jump time
  for the equilibrium process $\mathscr L^{\mathrm N}$. By the
  independence properties of the Poisson processes which generate
  $\mathscr T^{\mathrm N}$, we have $\mathscr T_{F-}^{\mathrm N}
  \stackrel{d}{=} \mathscr T_0^{\mathrm N}$ and consequently $\mathscr
  L_{F-}^{\mathrm N} \stackrel{d}{=} \mathscr L_0^{\mathrm N}
  $. Moreover, the jump removes a randomly chosen \emph{external
    branch} from $\mathscr T_{F-}^{\mathrm N}$. Since the $\mathrm
  N$-coalescent, restricted to $\mathrm N-1$ randomly chosen
  individuals, is in distribution identical to the $\mathrm
  (N-1)$-coalescent, it follows that $\mathscr T_{F}^{\mathrm N}
  \stackrel{d}{=} \mathscr T_0^{\mathrm N-1}$ and consequently
  $\mathscr L_{F}^{\mathrm N} \stackrel{d}{=} \mathscr L_0^{\mathrm
    N-1}+ 2\log\big( 1 - \tfrac 1{\mathrm N}\big)$. Moreover, the jump
  size, given by $\mathscr L_{F-}^{\mathrm N}-\mathscr L_{F}^{\mathrm
    N}$, is in distribution identical to a randomly chosen external
  branch of a $\mathrm N$-coalescent in equilibrium. Properties of the
  external branch length distribution are recalled in Section
  \ref{S:ext} and were studied in more detail in
  \cite{CaliebeNeiningerKrawczakRoesler2007}. For our setup, these
  results imply
$$ \mathrm N ( \mathscr L_{F-}^{\mathrm N}-\mathscr L_{F}^{\mathrm N}) 
\xRightarrow{\mathrm N\to\infty} J$$ for some random variable $J$,
taking values in the positive reals with expectation 2 and density
$x\mapsto 8/(2+x)^3$. This power law with exponent 3 was already
guessed in \cite{RauchBarYam2004} based on simulations. The random
variable $J$ has unbounded variance. For the asymptotics of this
variance, \cite{FuLi1993} already showed (see also Proposition
\ref{P:ext}) that\footnote{For sequences
  $(a_{\mathrm N})_{\mathrm N=1,2,...}$ and $(b_{\mathrm N})_{\mathrm
    N=1,2,...}$ we write $a_{\mathrm N} \stackrel{\mathrm N\to\infty}
  \sim b_{\mathrm N}$ iff $a_{\mathrm N} / b_{\mathrm N}
  \xrightarrow{\mathrm N\to\infty} 1$.}
\begin{align}\label{eq:J}
  \mathbb V[\mathrm N ( \mathscr L_{F-}^{\mathrm N}-\mathscr L_{F}^{\mathrm N})] 
  \stackrel{\mathrm N\to\infty}\sim 8\log(\mathrm N).
\end{align}
\end{remark}

\begin{remark}[Heuristics on jump sizes]
  The fact that the approximate size of a randomly chosen jump is of
  the order $2/\mathrm N$ can also be seen from the dynamics of
  $\mathscr L^{\mathrm N}$. In one time unit, tree length is gained by
  growth of the tree at constant speed $\mathrm N$. Moreover, the
  process $\mathscr T^{\mathrm N}$ makes approximately $\binom{\mathrm
    N} 2$ jumps. Since $\mathscr L^{\mathrm N}$ is in equilibrium, the
  tree growth and the jumps have to compensate each other. Therefore
  the expected size of a single jump must be $\frac 2{N-1}$; compare
  also with Proposition \ref{P:ext}.

  In the light of \eqref{eq:J}, the fact that the limit process
  $\mathscr L$ has infinite infinitesimal variance would not be
  surprising if there were no dependencies between jump sizes: In a
  short time $t$, the process $\mathscr L^{\mathrm N}$ makes
  approximately $\binom{\mathrm N}{2}t$ downward jumps as $\mathrm
  N\to\infty$. If jumps would be independent, we would get from
  \eqref{eq:J} that the variance of $\mathscr L^{\mathrm N}_t -
  \mathscr L^{\mathrm N}_0$ is approximately $4t \log\mathrm N$.
\end{remark}

\begin{figure}
\begin{center}
\includegraphics[width=12cm]{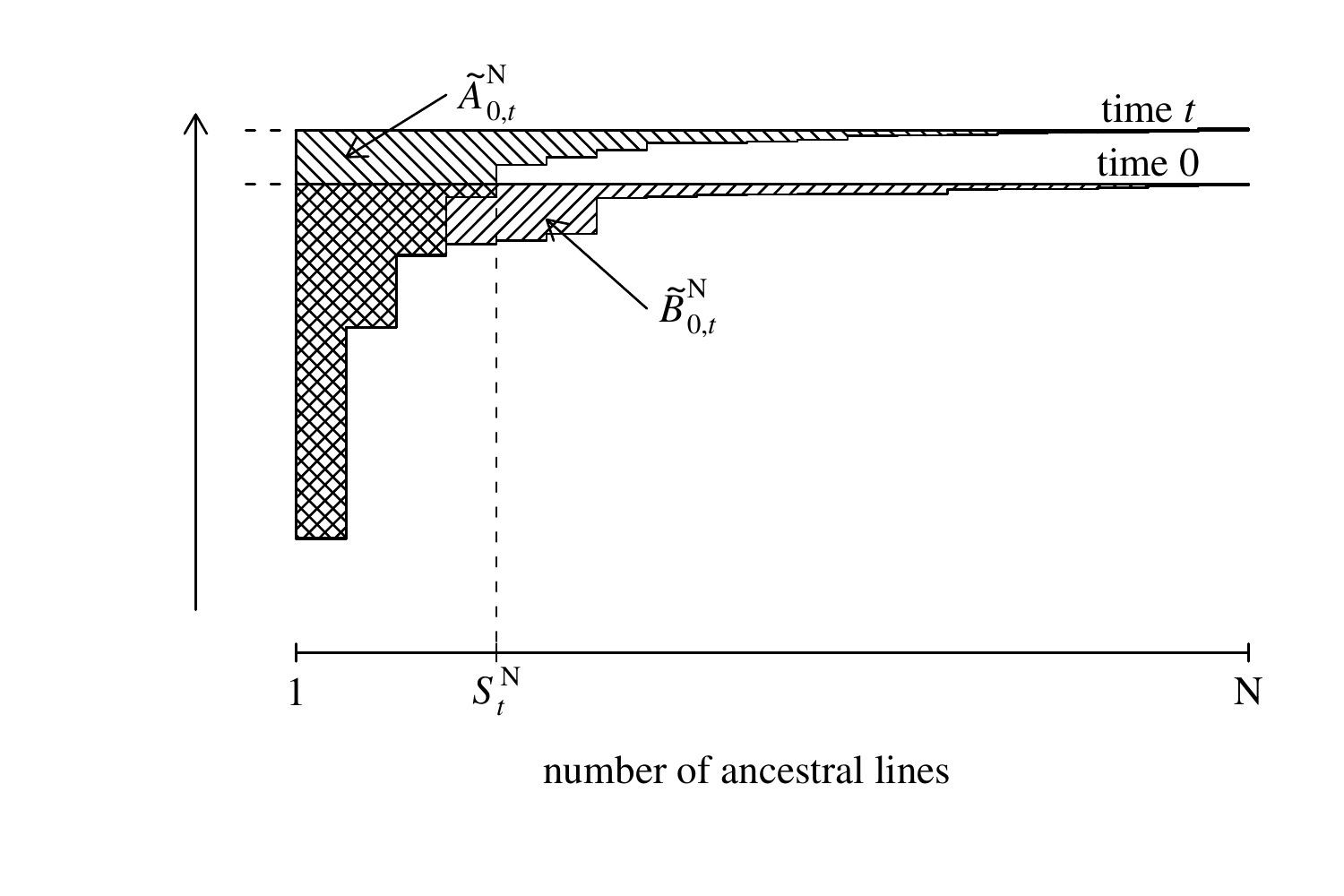}
\end{center}
\vspace{-1cm}

\caption{\label{fig1a} Schematic picture of tree change between two
  times 0 and $t$ for a population of size $\mathrm N$. Tree topology
  is ignored in the figure and only the number of ancestral lines is
  given. The population at time $t$ has $S_t^{\mathrm N}$ ancestors at
  time 0. The genealogical tree at time $t$ overlaps with the tree at
  time 0. The subset of the time-$t$ tree which does not belong to the
  time-$0$ tree is $\tilde A_{0,t}^{\mathrm N}$. Reversely, the subset
  of the time-$0$ tree which is lost between time 0 and $t$ is $\tilde
  B_{0,t}^{\mathrm N}$. The net difference in tree length is $\mathscr
  L_t^{\mathrm N}-\mathscr L_0^{\mathrm N} \stackrel{d}{=}\tilde
  A_{0,t}^{\mathrm N} - \tilde B_{0,t}^{\mathrm N}$.}
\end{figure}

\begin{remark}[Idea of the proof of  Theorem \ref{TH1}]
 \label{rem:idea}
 The crucial step in the proof, whose details are given in Section
 \ref{S:Sec4}, is to establish tightness of the family $\mathscr
 L^{\mathrm N}$. This uses auxiliary calculations on the evolution of
 Moran models (Section \ref{Moran}): it suffices to show that 'large'
 jumps in the tree length do not happen too often. To be more
 specific, we must show that during times $(t-h,t]$ and $(t,t+h]$ some
 moment of the smaller jump, $\mathscr L^{\mathrm N}_{t}-\mathscr
 L^{\mathrm N}_{t-h}$ or $\mathscr L^{\mathrm N}_{t+h}-\mathscr
 L^{\mathrm N}_{t}$, is bounded by $Ch^\theta$ for some constant $C$
 and $\theta>1$. Heuristically, such a statement is true since it can
 be shown that the times at which one of the $f$ oldest families of
 the coalescent tree dies out build a Poisson process with rate
 $\binom{f}{2}$, $f=2,3,\ldots$ (see Lemma \ref{L:moran}). Dying out
 of one of the $f$ oldest families implies larger jumps for smaller
 $f$ and the proof of tightness requires bounds for $f$ depending on
 the time interval $h$, using that loss times for one of the $f$
 oldest families in $(t-h,t]$ and $(t,t+h]$ are independent.

 To obtain the form for the infinitesimal variance, it is essential to
 bound jumps between times $0$ and $t$ of $\mathscr L^{\mathrm N}$ for
 small $t$, uniformly in $\mathrm N$; see Figure \ref{fig1a} for an
 illustration. Our proof is based on auxiliary calculations made in
 Sections \ref{sec:subtrees} and \ref{sec:ancestors}. Note that
 changes in $\mathscr L^{\mathrm N}$ come from two sources. First,
 between $0$ and $t$, additional tree length is gained by tree growth
 ($\tilde A_{0,t}^{\mathrm N}$ in the figure). The random variable
 $\tilde A_{0,t}^{\mathrm N}$ equals the tree-length of a Kingman
 $\mathrm N$-coalescent gained by time $t$ and for large $\mathrm N$,
 we see from Lemma \ref{L:treeLengthsSmallTimes} that $\lim_{\mathrm
   N\to\infty}\mathbb V[\tilde A_{0,t}^{\mathrm N}] \stackrel{t\to
   0}\sim \tfrac 23t$. Second, a part of the tree at time $t$ breaks
 off ($\tilde B_{0,t}^{\mathrm N}$ in the figure). This part is
 determined by the number $S_t^{\mathrm N}$ of ancestors at time 0 of
 the population at time $t$. Additionally, note that the number of
 ancestors at time $t$ of the population of size $\mathrm N$ converges
 in distribution to some random variable $S_t$ with $S_t
 \stackrel{t\to 0}\sim \lfloor \tfrac 2 t\rfloor$; see Lemma
 \ref{L:At}. Denoting by $B_{\mathrm K}$ the difference of the
 compensated tree length of a coalescent with infinitely many lineages
 and the compensated length of the tree spanned by a subset of
 $\mathrm K$ lineages, using some regularity, we can then show that
 $\mathbb V[B_{0,t}^{\mathrm N}] \stackrel{\mathrm N\to\infty}\sim
 \mathbb V[B_{S_t}] \stackrel{t\to 0}\sim \mathbb V[B_{\lfloor \tfrac
   2t\rfloor}] \stackrel{t\to 0} \sim 4t \log \tfrac 2t\stackrel{t\to
   0} \sim 4t |\log t|$ by Proposition \ref{P:34}.
 Combining the results for $\tilde A_{0,t}^{\mathrm N}$ and $\tilde
 B_{0,t}^{\mathrm N}$ we see that $\lim_{\mathrm N\to\infty}\mathbb
 V[\mathscr L_t^{\mathrm N} - \mathscr L_0^{\mathrm N}] \stackrel{t\to
   0} \sim 4 t | \log t|$.
\end{remark}


\begin{remark}[Convergence of tree-valued processes]
  In \cite{GrevenPfaffelhuberWinter2010}, a topology $\tau$ on the
  space of trees is specified and it is proved that the sequence of
  tree-valued processes $\mathscr T^{\mathrm N}$ converges in
  distribution to a tree-valued process $\mathscr T$, whose paths are
  a.s.\ continuous with respect to the topology $\tau$. One might be
  tempted to use convergence of $\mathscr T^{\mathrm N}$ to $\mathscr
  T$ in order to show that $\mathscr L^{\mathrm N}$ converges to some
  limit process $\mathscr L$ as well. Such an attempt would require
  that the function $\ell$ mapping finite trees to their (compensated)
  lengths is $\tau$-continuous. However, if $\ell$ would be
  $\tau$-continuous we would conclude that $\mathscr L$ has continuous
  paths, but $\mathscr L$ clearly makes jumps. Hence, $\ell$ is not
  continuous and convergence of $\mathscr T^{\mathrm N}$ cannot be
  used to show convergence of $\mathscr L^{\mathrm N}$.


\end{remark}

\begin{remark}[Extension to $\Lambda$-coalescents]
  In the past decade, the so-called $\Lambda$-coalescents
  \cite{Pitman1999} have gained increasing interest (see e.g.\
  \cite{Berest09} and references therein). Each of these coalescent
  processes arises as the large population limit of the genealogy of
  Cannings models and is uniquely determined by a finite measure
  $\Lambda$ on $[0;1]$. The Kingman coalescent then arises for
  $\Lambda=\delta_0$. For $\Lambda\neq \delta_0$, the underlying
  Cannings models have unbounded variance and the
  $\Lambda$-coalescents admits the possibility of more than two lines
  merging at the same time. An interesting direction for further
  research is the investigation of potential limits of the tree-length
  process of the respective genealogies of Cannings-models. For this,
  recent results by \cite{Moehle2006}, \cite{DrmotaEtAl2007},
  \cite{BBS2009} and \cite{DelmasEtAl2008} on moments and rescalings
  of $\Lambda$-coalescent trees, and by \cite{BBL2010} on the speed of
  coming down from infinity (which extends Aldous' result for the
  Kingman case, Lemma 4.6 below, to the $\Lambda$-case) will provide
  important ingredients.
\end{remark}

\begin{remark}[Connection to empirical population genetics]
  Coalescent trees are of particular importance in empirical
  population genetics and in the analysis of sequence diversity
  data. In the infinite sites model, mutations leading to segregating
  sites fall on the genealogical tree at constant rate. As a
  consequence, the number of segregating sites is Poisson distributed
  with a parameter proportional to the tree length. As illustrated by
  Figure \ref{fig:2}, the tree length process makes
  jumps. Particularly large jumps occur when the most recent common
  ancestor of the total population changes. At such a time $F$, one of
  the two oldest families in the population dies out and a long
  external branch breaks off the genealogical tree (see also
  \cite{Tajima1990}). At time $F-$ there are several segregating sites
  which are carried by all individuals which belong to the family
  which does not die out. Such segregating sites become fixed in the
  population when the MRCA of the population changes. In particular,
  fixation of segregating sites (also denoted by \emph{substitutions})
  come in bursts as time evolves, an observation already made by
  \cite{Watterson1982}. In addition, segregating sites which are
  present only in the oldest family which dies out at time $F$, are
  lost.

 Observations concerned with the substitutions of segregating sites
 are special properties of the \emph{mutation-drift balance}.  This
 dynamic equilibrium is between the introduction of new segregating
 sites due to mutation and loss of present ones due to genetic
 drift. Considered between times $0$ and $t$, the introduction of new
 mutations in the population are due to mutation events falling on the
 part of the genealogical tree gained between times $0$ and $t$ while
 the loss of existing mutations is due to some part of the
 genealogical tree at time $0$ breaking off by time $t$. Most
 interestingly, the number of segregating sites in the total
 population, unlike many other processes in population genetics, is
 'super-diffusive' in that it has infinite infinitesimal variance, as
 stated in Theorem \ref{TH1}.
\end{remark}

\section{Strong convergence of tree lengths}
\label{S:strong}
Proposition \ref{L:Gumbel} and Theorem \ref{TH1} establish convergence
in distribution for the real-valued random variables $\mathscr
L^{\mathrm N}_t$ ($t$ fixed) and the $\mathbb D$-valued random
variables $\mathscr L^{\mathrm N}$. We extend these results by
stronger notions of convergence, i.e.\ convergence in probability,
almost sure convergence and convergence in $L^2$. We start with
extensions for fixed times (Proposition \ref{L:Gumbel2}) and then come
to the extension involving the processes $\mathscr L^{\mathrm N}$
(Proposition \ref{P:TH1}). The proofs are given in Section
\ref{S:proofs2}.

\subsection{Tree lengths at fixed times}
For the extension of Proposition \ref{L:Gumbel}, fix $t\in\mathbb
R$. An elegant way to encode a random coalescent tree $\mathcal
T:=\mathscr T_t$ is in terms of (the completion of) a random metric on
$\mathbb N$, as proposed by Evans \cite{Evans2000}. To visualize this,
consider a sequence of lineages indexed by $\mathbb N$, where lineage
$i$ starts at time $0$ in leaf $i$. Any pair of lineages coalesces
independently at rate 1, and a random (ultra)metric $R$ is defined by
\begin{align} \label{eq:R} R(i,j) := 2\cdot \text{time to the most
    recent common ancestor of leaves }i,j.
\end{align}
The completion of $(\mathbb N,R)$ is a.s.\ a compact ultra-metric space
that represents the (uncountable set of) leaves of the coalescent tree. 

There are two canonical ways to approach the compensated length of
$\mathcal T$ by a sequence of lengths of finite trees.

\begin{figure}
\begin{center}
\includegraphics[width=13cm]{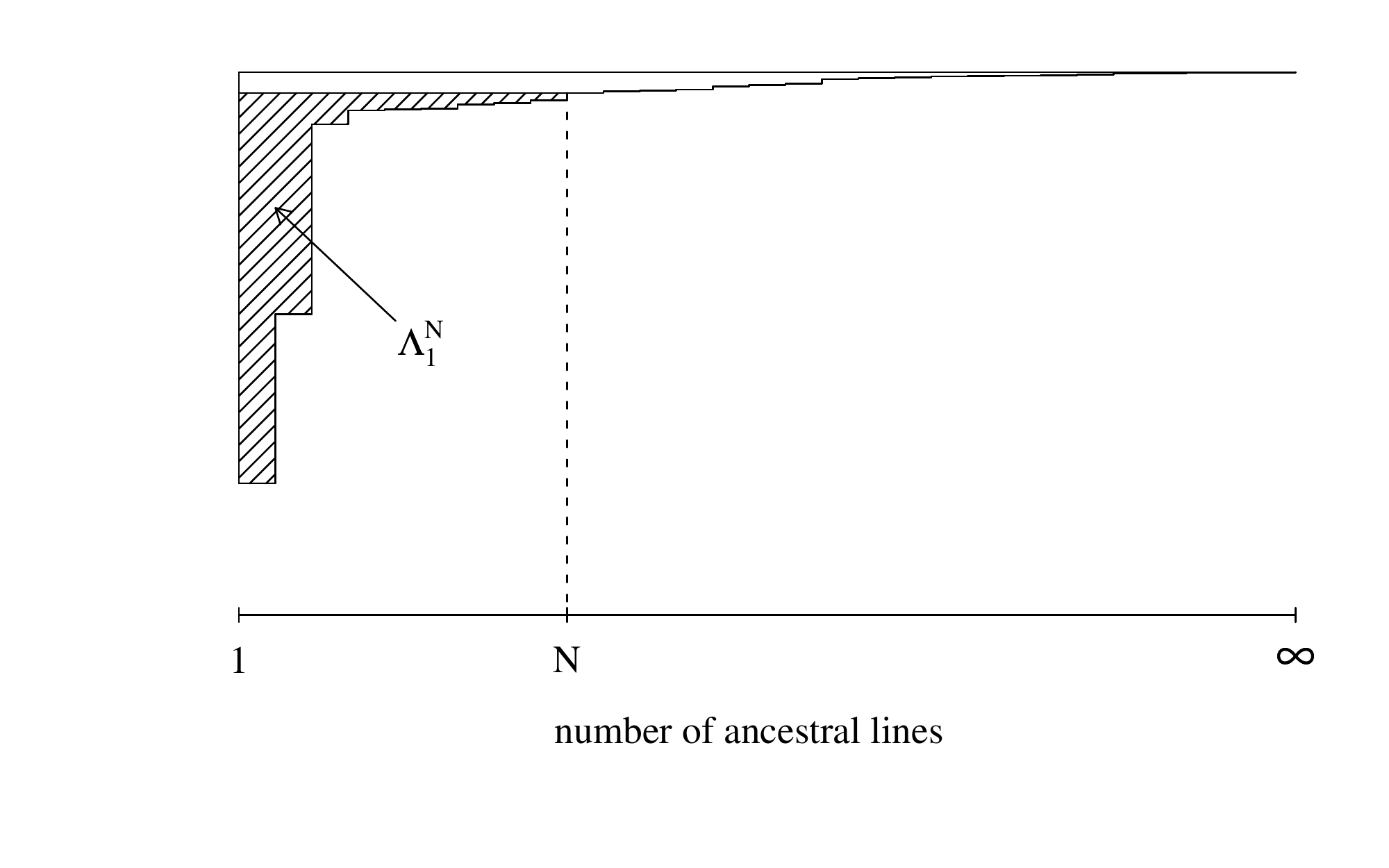}
\end{center}
\vspace{-1cm}

\caption{\label{fig:3} Kingman's temporal coupling: In the tree
  $\mathcal T$ coming down from infinitely many leaves, a tree with
  $\mathrm N$ leaves can be embedded by considering the part of
  $\mathcal T$ which is below the time at which $\mathcal T$ comes
  down to $\mathrm N$ lines.  The resulting tree length is
  $\Lambda_1^{\mathrm N}$. }
\end{figure}

\begin{figure}
\begin{center}
 \includegraphics[width=13cm]{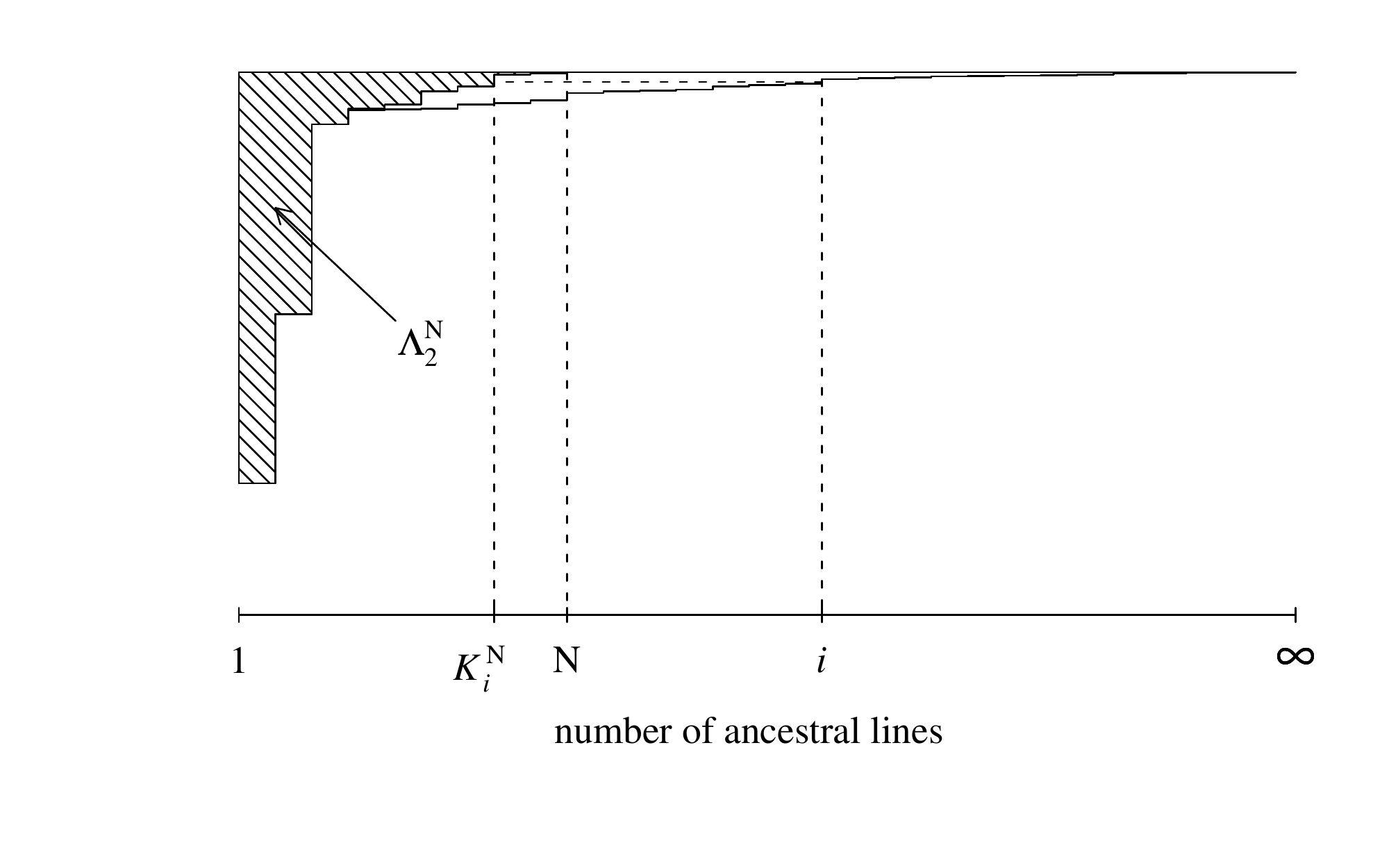}
\end{center}
\vspace{-1cm}

\caption{\label{fig:3b}Kingman's natural coupling: In the tree
  $\mathcal T$ coming down from infinitely many leaves, a tree with
  $\mathrm N$ leaves can be embedded by working down from the first
  $\mathrm N$ leaves (obtained by sampling) from the infinite
  tree. The resulting tree length is $\Lambda_2^{\mathrm N}$. In the
  natural coupling, $K_i^{\mathrm N}$ is the number of lines in the
  small tree while the full tree has $i$ lines.}
\end{figure}

\begin{enumerate}
\item \emph{From the root to the leaves:} Let $R_{(1)}, R_{(2)},
\ldots$ be a listing of the set $\{R(i,j): i, j \in \mathbb N, i\neq j\}$ in
decreasing order. The random variable
\begin{align}
  \label{eq:inter} X_k:= \tfrac 12(R_{(k-1)}-R_{(k)})
\end{align}
then gives the time the tree $\mathcal T$ spends with $k$ lines ``in
parallel'', $k=2,3,...$. We set
$$ \Lambda^{\mathrm N}_1 := \sum_{k=2}^{\mathrm N} kX_k - 2 \log \mathrm N, \quad {\mathrm N}=2,3,\ldots$$
The random sequence $(\Lambda_1^{\mathrm N})_{\mathrm N=2,3,...}$ is called
\emph{temporal coupling} in \cite{Kingman1982}; see also Figure~\ref{fig:3}.
\item \sloppy \emph{Across lineages:} For $\mathrm N=2,3,...$,
consider the finite subtree
\begin{align}\label{eq:defNat}
  \mathcal T^{\mathrm N}\text{ encoded by }\left(\{1,\ldots, N\},
    R\big |_{\{1,\ldots, N\}^2}\right), \qquad N=1,2\ldots.
\end{align}
Define
\begin{align*}
  \Lambda_2^{\mathrm N} := \widetilde \ell( \mathcal T^{\mathrm N})
  - 2 \log\mathrm N, \quad {\mathrm N}=2,3,\ldots
\end{align*}
This random sequence is called \emph{natural coupling} in
\cite{Kingman1982}; see also Figure \ref{fig:3b}.
\end{enumerate}

\begin{proposition}
\label{L:Gumbel2}
There is a random variable $\Lambda=\Lambda(R)$ such that $\tfrac 12
\Lambda$ is Gumbel distributed,
\begin{align}\label{eq:P11}
  \Lambda^{\mathrm N}_1 \xrightarrow{\mathrm N\to\infty} \Lambda
\end{align}
almost surely and in $L^2$ and
\begin{align}\label{eq:P12}
  \Lambda^{\mathrm N}_2 \xrightarrow{\mathrm N\to\infty} \Lambda
\end{align}
in $L^2$.
\end{proposition}

\subsection{The evolution of tree lengths}
\begin{figure}
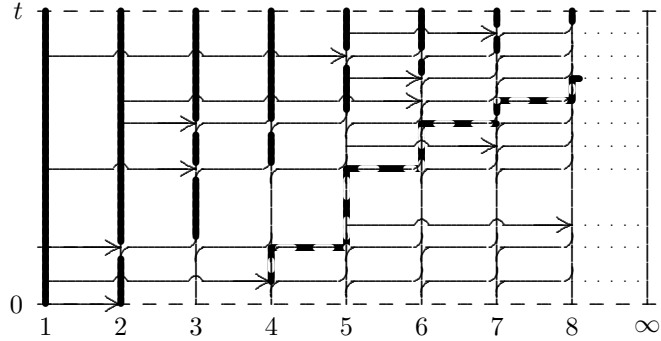

\label{lookdowngraph}
\begin{center}
\hspace{.5cm}
\beginpicture
\setcoordinatesystem units <1cm,1.5cm>
\setplotarea x from 1 to 9, y from 0.5 to 3.5
\multiput {{\color{black}{\tiny $\bullet$}}} at 1 3.4 *260 0 -.01 /
\multiput {{\color{black}{\tiny $\bullet$}}} at 2 3.4 *205 0 -.01 /
\multiput {{\color{black}{\tiny $\bullet$}}} at 2 1.2 *40 0 -.01 /
\multiput {{\color{black}{\tiny $\bullet$}}} at 3 3.4 *96 0 -0.01 /
\multiput {{\color{black}{\tiny $\bullet$}}} at 3 2.3 *25 0 -0.01 /
\multiput {{\color{black}{\tiny $\bullet$}}} at 3 1.9 *50 0 -0.01 /
\multiput {{\color{black}{\tiny $\bullet$}}} at 4 3.4 *96 0 -0.01 /
\multiput {{\color{black}{\tiny $\bullet$}}} at 4 2.3 *25 0 -0.01 /
\multiput {{\color{black}{\tiny $\bullet$}}} at 5 2.9 *38 0 -0.01 /
\multiput {{\color{black}{\tiny $\bullet$}}} at 5 3.4 *33 0 -0.01 /
\multiput {{\color{black}{\tiny $\bullet$}}} at 6 3.4 *33 0 -0.01 /
\multiput {{\color{black}{\tiny $\bullet$}}} at 6 2.95 *10 0 -0.01 /
\multiput {{\color{black}{\tiny $\bullet$}}} at 7 3.4 *12 0 -0.01 /
\multiput {{\color{black}{\tiny $\bullet$}}} at 7 3.15 *10 0 -0.01 /
\multiput {{\color{black}{\tiny $\bullet$}}} at 8 3.4 *10 0 -0.01 /

\multiput {{\color{black}{\tiny $\bullet$}}} at 4 1.0 *30 0 0.01 /
\multiput {{\color{white}{\tiny $\bullet$}}} at 4 1.1 *8 0 0.01 /
\multiput {{\color{black}{\tiny $\bullet$}}} at 4 1.3 *100 0.01 0 /
\multiput {{\color{white}{\tiny $\bullet$}}} at 4.1 1.3 *20 0.01 0 /
\multiput {{\color{white}{\tiny $\bullet$}}} at 4.6 1.3 *20 0.01 0 /

\multiput {{\color{black}{\tiny $\bullet$}}} at 5 1.3 *70 0 .01 /
\multiput {{\color{white}{\tiny $\bullet$}}} at 5 1.4 *15 0 .01 /
\multiput {{\color{white}{\tiny $\bullet$}}} at 5 1.75 *15 0 .01 /
\multiput {{\color{white}{\tiny $\bullet$}}} at 5 1.3 *1 0.01 0 /

\multiput {{\color{black}{\tiny $\bullet$}}} at 5 2.0 *100 0.01 0 /
\multiput {{\color{white}{\tiny $\bullet$}}} at 5.2 2.0 *20 0.01 0 /
\multiput {{\color{white}{\tiny $\bullet$}}} at 5.8 2.0 *10 0.01 0 /

\multiput {{\color{black}{\tiny $\bullet$}}} at 6 2.0 *40 0 .01 /
\multiput {{\color{white}{\tiny $\bullet$}}} at 6 2.2 *15 0 .01 /
\multiput {{\color{white}{\tiny $\bullet$}}} at 6 2 *5 0 0.01 /

\multiput {{\color{black}{\tiny $\bullet$}}} at 6 2.4 *100 0.01 0 /
\multiput {{\color{white}{\tiny $\bullet$}}} at 6.2 2.4 *20 0.01 0 /
\multiput {{\color{white}{\tiny $\bullet$}}} at 6.7 2.4 *20 0.01 0 /

\multiput {{\color{black}{\tiny $\bullet$}}} at 7 2.5 *10 0 .01 /

\multiput {{\color{black}{\tiny $\bullet$}}} at 7 2.6 *100 0.01 0 /
\multiput {{\color{white}{\tiny $\bullet$}}} at 7.2 2.6 *20 0.01 0 /
\multiput {{\color{white}{\tiny $\bullet$}}} at 7 2.6 *1 0.01 0 /
\multiput {{\color{white}{\tiny $\bullet$}}} at 7.7 2.6 *20 0.01 0 /

\multiput {{\color{black}{\tiny $\bullet$}}} at 8 2.6 *20 0 .01 /
\multiput {{\color{white}{\tiny $\bullet$}}} at 8 2.68 *5 0 .01 /

\multiput {{\color{black}{\tiny $\bullet$}}} at 8 2.8 *10 0.01 0 /\put{$1$} [cC] at 1 0.6
\put{$2$} [cC] at 2 0.6
\put{$3$} [cC] at 3 0.6
\put{$4$} [cC] at 4 0.6
\put{$5$} [cC] at 5 0.6
\put{$6$} [cC] at 6 0.6
\put{$7$} [cC] at 7 0.6
\put{$8$} [cC] at 8 0.6
\put{$\infty$} [cC] at 9 0.6

\setdashes
\plot 0.9 0.8 9.2 0.8 /
\plot 0.9 1.3 1.5 1.3 /
\plot 0.9 3.4 9.2 3.4 /

\put{$0$} [rC] at 0.7 0.8
\put{$t$} [rC] at 0.7 3.4

\setsolid
\thinlines

\plot 1 3.4 1 0.8 /
\plot 2 3.4 2 0.8 /
\plot 3 3.4 3 0.8 /
\plot 4 3.4 4 0.8 /
\plot 5 3.4 5 0.8 /
\plot 6 3.4 6 0.8 /
\plot 7 3.4 7 0.8 /
\plot 8 3.4 8 0.8 /
\plot 9 3.4 9 0.8 /

\setdots
\plot 8 1.0 9 1.0 /
\plot 8 1.3 9 1.3 /
\plot 8 1.5 9 1.5 /
\plot 8 2 9 2 /
\plot 8 2.2 9 2.2 /
\plot 8 2.4 9 2.4 /
\plot 8 2.6 9 2.6 /
\plot 8 2.8 9 2.8 /
\plot 8 3 9 3 /
\plot 8 3.2 9 3.2 /


\setsolid

\plot 1 1.0 1.9 1.0 /
\plot 2.1 1.0 2.9 1.0 /
\plot 3.1 1.0 3.9 1.0 /
\plot 4.1 1.0 4.9 1.0 /
\plot 5.1 1.0 5.9 1.0 /
\plot 6.1 1.0 6.9 1.0 /
\plot 7.1 1.0 7.9 1.0 /

\plot 1 1.3 1.9 1.3 /
\plot 2.1 1.3 2.9 1.3 /
\plot 3.1 1.3 3.9 1.3 /
\plot 4.1 1.3 4.9 1.3 /
\plot 5.1 1.3 5.9 1.3 /
\plot 6.1 1.3 6.9 1.3 /
\plot 7.1 1.3 7.9 1.3 /

\plot 5 1.5 5.9 1.5 /
\plot 6.1 1.5 6.9 1.5 /
\plot 7.1 1.5 7.9 1.5 /


\plot 1 2 1.9 2 /
\plot 2.1 2 2.9 2 /
\plot 3.1 2 3.9 2 /
\plot 4.1 2 4.9 2 /
\plot 5.1 2 5.9 2 /
\plot 6.1 2 6.9 2 /
\plot 7.1 2 7.9 2 /

\plot 5 2.2 5.9 2.2 /
\plot 6.1 2.2 6.9 2.2 /
\plot 7.1 2.2 7.9 2.2 /

\plot 2 2.4 2.9 2.4 /
\plot 3.1 2.4 3.9 2.4 /
\plot 4.1 2.4 4.9 2.4 /
\plot 5.1 2.4 5.9 2.4 /
\plot 6.1 2.4 6.9 2.4 /
\plot 7.1 2.4 7.9 2.4 /

\plot 2 2.6 2.9 2.6 /
\plot 3.1 2.6 3.9 2.6 /
\plot 4.1 2.6 4.9 2.6 /
\plot 5.1 2.6 5.9 2.6 /
\plot 6.1 2.6 6.9 2.6 /
\plot 7.1 2.6 7.9 2.6 /

\plot 5 2.8 5.9 2.8 /
\plot 6.1 2.8 6.9 2.8 /
\plot 7.1 2.8 7.9 2.8 /

\plot 1 3 1.9 3 /
\plot 2.1 3 2.9 3 /
\plot 3.1 3 3.9 3 /
\plot 4.1 3 4.9 3 /
\plot 5.1 3 5.9 3 /
\plot 6.1 3 6.9 3 /
\plot 7.1 3 7.9 3 /


\plot 5 3.2 5.9 3.2 /
\plot 6.1 3.2 6.9 3.2 /
\plot 7.1 3.2 7.9 3.2 /

\arrow <0.2cm> [0.375,1] from 6.5 3.2 to 7 3.2
\arrow <0.2cm> [0.375,1] from 4.5 3 to 5 3
\arrow <0.2cm> [0.375,1] from 5.5 2.8 to 6 2.8
\arrow <0.2cm> [0.375,1] from 5.5 2.6 to 6 2.6
\arrow <0.2cm> [0.375,1] from 2.5 2.4 to 3 2.4
\arrow <0.2cm> [0.375,1] from 6.5 2.2 to 7 2.2
\arrow <0.2cm> [0.375,1] from 2.5 2 to 3 2
\arrow <0.2cm> [0.375,1] from 7.5 1.5 to 8 1.5
\arrow <0.2cm> [0.375,1] from 1.5 1.3 to 2 1.3
\arrow <0.2cm> [0.375,1] from 3.5 1 to 4 1
\arrow <0.2cm> [0.375,1] from 1 0.8 to 2 0.8

\circulararc 90 degrees from 2.1 1.3 center at 2.15 1.2

\circulararc 90 degrees from 2.9 1.3 center at 2.85 1.4

\circulararc 90 degrees from 3.1 2 center at 3.15 1.9
\circulararc 90 degrees from 3.1 2.4 center at 3.15 2.3
\circulararc 90 degrees from 3.1 1.3 center at 3.15 1.2

\circulararc 90 degrees from 3.9 2 center at 3.85 2.1
\circulararc 90 degrees from 3.9 2.4 center at 3.85 2.5
\circulararc 90 degrees from 3.9 1.3 center at 3.85 1.4

\circulararc 90 degrees from 4.1 1 center at 4.15 0.9
\circulararc 90 degrees from 4.1 1.3 center at 4.15 1.2
\circulararc 90 degrees from 4.1 2 center at 4.15 1.9
\circulararc 90 degrees from 4.1 2.4 center at 4.15 2.3

\circulararc 90 degrees from 4.9 1.0 center at 4.85 1.1
\circulararc 90 degrees from 4.9 1.3 center at 4.85 1.4
\circulararc 90 degrees from 4.9 2 center at 4.85 2.1
\circulararc 90 degrees from 4.9 2.4 center at 4.85 2.5

\circulararc 90 degrees from 5.1 1.0 center at 5.15 0.9
\circulararc 90 degrees from 5.1 1.3 center at 5.15 1.2
\circulararc 90 degrees from 5.1 2 center at 5.15 1.9
\circulararc 90 degrees from 5.1 2.4 center at 5.15 2.3
\circulararc 90 degrees from 5.1 3.0 center at 5.15 2.9

\circulararc 90 degrees from 5.9 1.0 center at 5.85 1.1
\circulararc 90 degrees from 5.9 1.3 center at 5.85 1.4
\circulararc 90 degrees from 5.9 2 center at 5.85 2.1
\circulararc 90 degrees from 5.9 2.4 center at 5.85 2.5
\circulararc 90 degrees from 5.9 3.0 center at 5.85 3.1

\circulararc 90 degrees from 6.1 1.0 center at 6.15 0.9
\circulararc 90 degrees from 6.1 1.3 center at 6.15 1.2
\circulararc 90 degrees from 6.1 2 center at 6.15 1.9
\circulararc 90 degrees from 6.1 2.4 center at 6.15 2.3
\circulararc 90 degrees from 6.1 2.6 center at 6.15 2.5
\circulararc 90 degrees from 6.1 2.8 center at 6.15 2.7
\circulararc 90 degrees from 6.1 3.0 center at 6.15 2.9

\circulararc 90 degrees from 6.9 1.0 center at 6.85 1.1
\circulararc 90 degrees from 6.9 1.3 center at 6.85 1.4
\circulararc 90 degrees from 6.9 2 center at 6.85 2.1
\circulararc 90 degrees from 6.9 2.4 center at 6.85 2.5
\circulararc 90 degrees from 6.9 2.6 center at 6.85 2.7
\circulararc 90 degrees from 6.9 2.8 center at 6.85 2.9
\circulararc 90 degrees from 6.9 3.0 center at 6.85 3.1

\circulararc 90 degrees from 7.1 1.0 center at 7.15 0.9
\circulararc 90 degrees from 7.1 1.3 center at 7.15 1.2
\circulararc 90 degrees from 7.1 2 center at 7.15 1.9
\circulararc 90 degrees from 7.1 2.2 center at 7.15 2.1
\circulararc 90 degrees from 7.1 2.4 center at 7.15 2.3
\circulararc 90 degrees from 7.1 2.6 center at 7.15 2.5
\circulararc 90 degrees from 7.1 2.8 center at 7.15 2.7
\circulararc 90 degrees from 7.1 3.0 center at 7.15 2.9
\circulararc 90 degrees from 7.1 3.2 center at 7.15 3.1

\circulararc 90 degrees from 7.9 1.0 center at 7.85 1.1
\circulararc 90 degrees from 7.9 1.3 center at 7.85 1.4
\circulararc 90 degrees from 7.9 2 center at 7.85 2.1
\circulararc 90 degrees from 7.9 2.2 center at 7.85 2.3
\circulararc 90 degrees from 7.9 2.4 center at 7.85 2.5
\circulararc 90 degrees from 7.9 2.6 center at 7.85 2.7
\circulararc 90 degrees from 7.9 2.8 center at 7.85 2.9
\circulararc 90 degrees from 7.9 3.0 center at 7.85 3.1
\circulararc 90 degrees from 7.9 3.2 center at 7.85 3.3

\setquadratic



\plot 1.9 1.0 2 1.05 2.1 1.0 /
\plot 2.9 1.0 3 1.05 3.1 1.0 /

\plot 5.9 1.5 6 1.55 6.1 1.5 /
\plot 6.9 1.5 7 1.55 7.1 1.5 /


\plot 1.9 2 2 2.05 2.1 2 /

\plot 5.9 2.2 6 2.25 6.1 2.2 /

\plot 2.9 2.6 3 2.65 3.1 2.6 /
\plot 3.9 2.6 4 2.65 4.1 2.6 /
\plot 4.9 2.6 5 2.65 5.1 2.6 /

\plot 1.9 3 2 3.05 2.1 3 /
\plot 2.9 3 3 3.05 3.1 3 /
\plot 3.9 3 4 3.05 4.1 3 /


\plot 5.9 3.2 6 3.25 6.1 3.2 /

\endpicture
\end{center}
\caption{\label{FigLookdown0}{Detail of a look-down graph. Time is
  running upwards; all lines at the first 8 levels are drawn between
  times $0$ and $t$. At times in $\mathcal P_{ij}$ an arrow is drawn
  from $i$ to $j$, and all lines at levels  $\ge j$ are pushed
  upwards as indicated by bent lines. The dashed line initiates one
  partition element of the partition $\mathcal G$ of $\mathbb
  R\times \{2,3,\ldots\}$. The solid lines constitute
  the tree $ \mathcal T^{\mathrm {ld},\mathrm 8}_t$ of the first 8 ancestral lineages back from time $t$. The sum of the lengths of the 
  solid lines is the tree length of the first $8$ levels at time
  $t$.}}
\end{figure}

For the extension of Theorem \ref{TH1} we briefly review the lookdown
construction of \cite{DonnellyKurtz1999}; see also
\cite{PfaffelhuberWakolbinger2006} for a detailed description and
Figure \ref{FigLookdown0} for an illustration.  

Consider the set of vertices $\mathcal V := \mathbb R\times\mathbb
N$. A vertex $(t,i)$ is referred to as {\em the individual at time $t$
  at level $i$}. The source of randomness in the lookdown construction
is a family of rate one Poisson processes $(\mathcal P_{ij})_{1\leq
  i<j}$. At times $t\in\mathcal P_{ij}$, the individual at level $j$
\emph{looks down} to level $i$. As the illustration in Figure
\ref{FigLookdown0} shows, at a lookdown event in $\mathcal P_{ij}$,
all individuals at levels $k\geq j$ are pushed one level up, and a new
{\em line of ascent} is born at level $j$. Note that the individual at
level $k$ is pushed to level $k+1$ at rate $\binom k 2$.

We define the partition $\mathcal G$ of $\mathbb R \times
\{2,3,\ldots\}$ into lines of ascent as follows. Every $s_0\in\mathcal
P_{ij}$ gives rise to a partition element $G$ of the form
$$ ([s_0, s_1)\times \{j\} \cup 
([s_1,s_2)\times \{j+1\}) \cup ([s_2,s_3)\times \{j+2\}) \cup ...$$
with $s_{k+1}>s_k$ for all $k=0,1,2,...$. Here, $s_{k+1}$ is the
smallest element in $[s_k,\infty) \cap \bigcup_{1\leq i'<j'\leq
 k+j}\mathcal P_{i'j'}$. We say that $G$ is born by $(s_0,i)$ and
pushed one level up at times $s_1,s_2,...$. If $(s_0,i)\in G'$, we say
that $G$ descends from $G'$. Since the individual at level $k$ is
pushed up by one at rate $\binom k 2$ and hence, for $G$ as above,
$\lim_{k\to\infty}s_k$ is finite.

For $s\leq t,$ $i\leq j$ let $G, G'\in\mathcal G$ be such that
$(s,i)\in G', (t,j)\in G$. We say that $(s,i)$ is ancestor of $(t,j)$
if either $G=G'$ or there are $G_1,...,G_n$ such that $G$ descends
from $G_n$, $G_k$ descends from $G_{k-1}$, $k=2,...,n$ and $G_1$
descends from $G'$. In this case, we define $A_s(t,j):=i$. In
addition, for $t\in\mathbb R$, define the random metric $R_t^{\mathrm
  {ld}}$ (compare with \eqref{eq:R}) on $\mathbb N$ by
$$ R^{\mathrm {ld}}_t(i,j) = 2\cdot \inf\{t-s: A_s(t,i) = A_s(t,j)\}.$$
We define the random trees
$$ \mathcal T^{\mathrm {ld},\mathrm N}_t := \Big(\{1,...,\mathrm N\}, R_t^{\mathrm {ld}}
|_{\{1,...,\mathrm N\}^2}\Big), \qquad \mathrm N=2,3,...$$ and
$$ \mathscr L_t^{\mathrm {ld},\mathrm N} := \widetilde \ell(\mathcal T^{\mathrm {ld},\mathrm N}_t) - 2\log\mathrm N.$$
Recall the complete Skorokhod metric $d_{\mathrm {Sk}}$ on $\mathbb D$
from \cite[Section 3.5]{EthierKurtz1986}. We are now ready to state a
result extending Theorem \ref{TH1} to convergence in probability,
proved in Section \ref{S:proofs2}.

\begin{proposition}
\label{P:TH1}
There is a process $\mathscr L^{\mathrm {ld}}$, having the same distribution as $\mathscr L$
from Theorem \ref{TH1}, such that
$$ d_{\mathrm {Sk}}(\mathscr L^{\mathrm N,\mathrm {ld}}, \mathscr L^{\mathrm {ld}}) \xrightarrow{\mathrm N\to\infty} 0 $$
in probability.
\end{proposition}

\section[Kingman's coalescent and the Moran model]{Auxiliary results
  on Kingman's coalescent and the Moran model}
\label{S:3}
In this section we collect some facts on Kingman's coalescent
(Subsections \ref{S:41}--\ref{S:44}) and the Moran model (Subsection
\ref{Moran}) which will be required for the proof of Theorem
\ref{TH1}.

For the Kingman coalescent, we take the tree $\mathcal T$ as
introduced in Section \ref{S:strong}. Recall the subtrees $\mathcal
T^{\mathrm N}, \mathrm N=2,3,...$ as defined in \eqref{eq:defNat} and
the inter-coalescence times $X_2, X_3,...$ from
\eqref{eq:inter}. Recall that $X_k$ is exponentially distributed with
parameter $\binom k 2 $.

\subsection{The Markov Chain $K^{\mathrm N}$}
\label{S:41}
We define
\begin{align}
K^{\mathrm N}_i := \text{ number of lines in $\mathcal T^{\mathrm
    N}$ while $\mathcal T$ has $i$ lines}
\end{align}
and $K^{\mathrm N}:=(K_i^{\mathrm N})_{i=1,2,...}$; see also Figure
\ref{fig:3b}. Note that
\begin{align}
  \label{eq:lambda2}
  \Lambda_2^{\mathrm N} = \sum_{i=2}^\infty K_i^{\mathrm N}X_i - 2
  \log\mathrm N.
\end{align}
The connection between the trees $\mathcal T^{\mathrm N}$ and
$\mathcal T$ has been described e.g. by \cite{Tavare1984},
\cite{Saundersetal1984} and \cite[Section
4.4]{EtheridgePfaffelhuberWakolbinger2006}.  Lemma 4.8 of
\cite{EtheridgePfaffelhuberWakolbinger2006} states that $K^{\mathrm
  N}$ builds a Markov chain with one- and two-dimensional
distributions
\begin{align}\label{eq:distKi1}
 \mathbb P[K_i^{\mathrm N}=k] & = \frac{\binom{{\mathrm
       N}-1}{{\mathrm N}-k}\binom{i}{k}}{\binom{{\mathrm
       N}+i-1}{{\mathrm N}}}, \qquad i\geq 1, 1\leq k\leq \mathrm
 N\\\label{eq:distKi2} \mathbb P[K_j^{\mathrm N}=\ell|K^{\mathrm
   N}_i=k] & = \frac{\binom{{\mathrm N}-k}{{\mathrm
       N}-\ell}\binom{j+k-1}{i+l-1}}{\binom{{\mathrm N}+j-1}{{\mathrm
       N}+i-1}},\qquad 1\leq i\leq j, 1\leq k\leq \ell\leq \mathrm N.
\end{align}
We will need some moment properties of this Markov chain.  We leave
out the straightforward details of the proof.

\begin{lemma}\label{L1} For $i\geq 1$,
  \begin{align}
  \label{eq:L11a}
  \mathbb E[i-K_i^{\mathrm N}] & = \frac{i(i-1)}{\mathrm N+i-1} ,\\
  \intertext{and for $1\leq i\leq j$}
  \label{eq:L11b}
  \mathbb E[(i-K_i^{\mathrm N})(j-K_j^{\mathrm N})] & =
  \frac{i(i-1)j(j-1)}{(\mathrm N+i-1)(\mathrm N+j-1)} +
  \frac{i(i-1)\mathrm N(\mathrm N-1)}{(\mathrm N+j-1)(\mathrm
    N+i-1)(\mathrm N+i-2)}.
  \end{align}
\end{lemma}

\subsection{The length of an external branch}
\label{S:ext}
We will recall several facts of the length of a randomly chosen
external branch in an $\mathrm N$-coalescent $\mathcal T^{\mathrm N}$.
In this setting we take the inter-coalescence times $X_2^{\mathrm
  N},...,X_{\mathrm N}^{\mathrm N}$ such that $X_i^{\mathrm N}$ is
exponentially distributed with rate $\binom i 2$, $i=2,...,\mathrm
N$. We denote by $J^{\mathrm N}$ the length of a randomly chosen
external branch. The results we describe in this section are collected
from \cite{FuLi1993}, \cite{Durrett2008}, \cite{BlumFrancois2005} and
\cite{CaliebeNeiningerKrawczakRoesler2007} and stated here for
completeness.

We define
$$ F^{\mathrm N} := f \qquad \text{iff}\qquad J^{\mathrm N} = \sum_{k=f+1}^{\mathrm N} X_k, $$
i.e.\ $F^{\mathrm N}$ denotes the number of lines extant in the $\mathrm
N$-coalescent at the time at which the external branch connects to the tree. We
give a basic fact about $F^{\mathrm N}$ and properties of $J^{\mathrm
 N}$.

\begin{proposition}[External branches]
 \label{P:ext}
 For $f=1,...,\mathrm N-1$,
 \begin{align*}
   \mathbb P[F^{\mathrm N}<f] = \frac{f(f-1)}{\mathrm N(\mathrm N-1)}
   \qquad \text{i.e.} \qquad \mathbb P[F^{\mathrm N}=f] =
   \frac{2f}{\mathrm N(\mathrm N-1)}
 \end{align*}
 The first two moments of $J^{\mathrm N}$ are given by
 \begin{align*}
   \mathbb E[J^{\mathrm N}] & = \frac{2}{\mathrm N},\\
   \mathbb V[J^{\mathrm N}] & = \frac{8\sum_{k=1}^{\mathrm N}\tfrac
     1k - 12 + \tfrac{4}{\mathrm N}}{\mathrm N(\mathrm N-1)}.
 \end{align*}  
\end{proposition}

\begin{proof}
  We obtain the distribution of $F^{\mathrm N}$ as follows: With
  probability $\frac{\binom{\mathrm N-1}{2}}{\binom{\mathrm N}{2}}$
  the randomly chosen external branch is not involved in the first
  coalescence event (bringing the number of lines from $\mathrm N$
  down to $\mathrm N -1$). Iterating this argument, we immediately see
  that the probability that the randomly chosen line did not take part
  in the first $\mathrm N-f+1$ coalescence events is
 \begin{align*}
   \mathbb P[F^{\mathrm N}<f] = \frac{\binom{\mathrm
       N-1}{2}}{\binom{\mathrm N}{2}} \cdot \frac{\binom{\mathrm
       N-2}{2}}{\binom{\mathrm N-1}{2}} \cdots \frac{\binom{
       f}{2}}{\binom{f+1}{2}} = \frac{\binom{f}{2}}{\binom{\mathrm
       N}{2}} = \frac{f(f-1)}{\mathrm N(\mathrm N-1)}.
 \end{align*}
 To compute moments of $J^{\mathrm N}$, we use the representation
 \begin{align*}
   J^{\mathrm N} = \sum_{k=F^{\mathrm N}+1}^{\mathrm N} X_k^{\mathrm
     N}.
 \end{align*}
 Recalling that $X_2^{\mathrm N},...,X_{\mathrm N}^{\mathrm N}$ are
 independent of the tree topology, and $F^{\mathrm N}$ is measurable
 with respect to the tree topology, we get that $F^{\mathrm N},
 X_2^{\mathrm N},...,X_{\mathrm N}^{\mathrm N}$ are independent. The
 first two moments of $J^{\mathrm N}$ are now obtained by
 \begin{align*}
   \mathbb E[J^{\mathrm N}] & = \mathbb E\big[ \mathbb E[J^{\mathrm
     N}|F^{\mathrm N}]\big] = \mathbb E\Big[ \frac{2}{F^{\mathrm N}} -
   \frac{2}{\mathrm N}\Big] = \Big( 2\sum_{f=1}^{\mathrm N-1}
   \frac{f}{\binom{\mathrm N}{2}} \frac 1 f\Big) - \frac{2}{\mathrm N}
   = \frac{2}{\mathrm N},\\
   \intertext{which also implies that $\mathbb E[\tfrac{2}{F^{\mathrm
         N}}] = \tfrac{4}{\mathrm N}$ and} \mathbb V[J^{\mathrm N}] &
   = \mathbb E\big[ \mathbb V[J^{\mathrm N}|F^{\mathrm N}]\big] +
   \mathbb V\big[ \mathbb E[J^{\mathrm N}|F^{\mathrm N}]\big] =
   \mathbb E\Big[ \sum_{k=F^{\mathrm N}+1}^{\mathrm
     N}\frac{1}{\binom{k}{2}^2}\Big] + \mathbb
   V\Big[\frac{2}{F^{\mathrm N}}\Big] \\ & = \sum_{f=1}^{\mathrm
     N-1}\sum_{k=f+1}^{\mathrm N} \frac{2f}{\mathrm N(\mathrm N-1)}
   \frac{1}{\binom{k}{2}^2} + 4\sum_{f=1}^{\mathrm N-1}
   \frac{2f}{\mathrm N(\mathrm N-1)}\frac{1}{f^2} - \frac{16}{\mathrm
     N^2}\\ & = \frac{1}{\binom{\mathrm N}{2}} \sum_{k=2}^{\mathrm
     N}\sum_{f=1}^{k-1} f \frac{1}{\binom{k}{2}^2} +
   \frac{4}{\binom{\mathrm N}{2}}\sum_{f=1}^{\mathrm N-1}
   f\frac{1}{f^2} - \frac{16}{\mathrm N^2}\\ & = \frac{4}{\mathrm N^2}
   + \frac{8}{\mathrm N(\mathrm N-1)} \Big(\sum_{f=1}^{\mathrm N-1}
   \frac{1}{f}\Big) - \frac{16}{\mathrm N^2}
   \\
   & = \frac{8\sum_{k=1}^{\mathrm N}\tfrac 1k - 12 + \tfrac{4}{\mathrm
       N}}{\mathrm N(\mathrm N-1)}.
 \end{align*}
\end{proof}

\subsection{Subtrees of coalescents and their lengths}\label{sec:subtrees}
The aim of this section is to analyze the difference of the tree
lengths of $\mathcal T$ and of $\mathcal T^{\mathrm N}$. Since this
difference is infinite, we have to carry out a limiting procedure, compensating by the mean. Recall that the inter-coalescence times $X_2, X_3,...$ are
independent of the tree topology of $\mathcal T$ in general and of
$K_2^{\mathrm N}, K_3^{\mathrm N},...$ in particular. We study the
random variable
\begin{align}
\label{eq:BKN0}
B_{\mathrm N} := \sum_{i=2}^{\infty} \Big( (i-K_i^{\mathrm N}) X_i - 
\mathbb E[i-K_i^{\mathrm N}]\cdot\mathbb E[X_i]\Big) 
\end{align}
which is the compensated difference of the tree lengths of $\mathcal T$ and
$\mathcal T^{\mathrm N}$. 

\begin{remark}
  Let us first make sure that the infinite sum in the definition of
  $B_{\mathrm N}$ exists and has expectation zero. To see this, we fix
  $\mathrm M \in \{2,3,...\}$ and consider the sequences $(B_{\mathrm
    N, \mathrm M})_{\mathrm M=2,3,...}$ where $B_{\mathrm N, \mathrm
    M}$ is defined as $B_{\mathrm N}$ but with the sum ranging from
  $i=2$ to $i= \mathrm M$. For $\mathrm M<\mathrm M'$, by Lemma
  \ref{L1}, and using that $\mathbb E[X_iX_j]\leq 2\cdot \mathbb
  E[X_i] \cdot \mathbb E[X_j]$ for all $i,j = 2,3,...$
\begin{align*}
  \mathbb E & \big[\big(B_{\mathrm N, \mathrm M'} - B_{\mathrm N,
    \mathrm M}\big)^2\big] \\ & = \mathbb E\Big[ \Big(\sum_{i=\mathrm
    M+1}^{\mathrm M'} \Big((i-K^{\mathrm N}_i)X_i - \frac{2}{\mathrm
    N+i-1}\Big)\Big)^2\Big]  \\
  & \leq 8 \sum_{i=\mathrm M+1}^{\mathrm M'} \sum_{j=i}^{\mathrm M'}
  \mathbb E\Big[\Big( \frac{i-K_i^{\mathrm N}}{i(i-1)} -
  \frac{1}{\mathrm N+i-1}\Big) \Big( \frac{j-K_j^{\mathrm N}}{j(j-1)}
  - \frac{1}{\mathrm N+j-1}\Big)\Big] \\ & \leq 8 \sum_{i=\mathrm
    M+1}^{\mathrm M'} \sum_{j=i}^{\mathrm M'} \frac{1}{j(j-1)(\mathrm
    N+i-1)} \leq 8 \sum_{i=\mathrm M+1}^{\mathrm M'} \frac{1}{(i-1)^2}
  \\ & = \mathcal O(\tfrac{1}{\mathrm M}) \qquad \text{as $\mathrm M
    \to \infty$.}
\end{align*}
Hence, the sequence $(B_{\mathrm N, \mathrm M})_{\mathrm M=2,3,...}$
is Cauchy in $L^2$ and thus converges in $L^2$ to the limit variable
$B_{\mathrm N}$ defined in \eqref{eq:BKN0}. Furthermore we obtain that
${\mathbb E}[B_{N}]=0$ by continuity of the linear functional
${\mathbb E}$ on $L^2$.
\end{remark}

\begin{proposition}[Variance of the difference in length of an
  infinite and a finite coalescent in the natural
  coupling]\label{P:34}
 \begin{align}\label{eq:toshow1}
   \mathbb E[B_{\mathrm N}^2] & \stackrel{\mathrm N\to\infty}\sim
   \frac{8\log\mathrm N}{\mathrm N}.
 \end{align}
\end{proposition}

\begin{proof}
  First,
  $$ \mathbb E[X_iX_j] = \mathbb E[X_i]\cdot\mathbb E[X_j]\cdot (1+\delta_{ij}).$$
  In order to obtain \eqref{eq:toshow1}, by a straightforward
  calculation using Lemma \ref{L1},
  \begin{equation}
    \label{eq:toshow11}
    \begin{aligned}
      \mathbb E[B^2_{\mathrm N}]
      & = 8 \sum_{j=3}^\infty \sum_{i=2}^{j-1} \frac{\mathrm N(\mathrm
        N-1)}{j(j-1)(\mathrm N+j-1)(\mathrm N+i-1)(\mathrm N+i-2)} \\
      & \qquad\qquad\qquad + \sum_{i=2}^\infty \Big(\frac{8\cdot
        \mathrm N(\mathrm N-1)}{i(i-1)(\mathrm N+i-1)^2 (\mathrm
        N+i-2)} + \frac{4}{(\mathrm N+i-1)^2} \Big).
  \end{aligned}
\end{equation}
It is clear that the expression in the last line tends to 0 as
$C/\mathrm N$, for some $C>0$, as $\mathrm N\to\infty$. For the
expression in the next to last line, we obtain
\begin{align*}
  8\sum_{j=3}^\infty & \frac{\mathrm N(\mathrm N-1)}{j(j-1)(\mathrm
    N+j-1)} \sum_{i=2}^{j-1} \Big(\frac{1}{\mathrm N+i-2} -
  \frac{1}{\mathrm N+i-1}\Big) \\ 
  & = 8\sum_{j=3}^\infty \frac{(\mathrm N-1)(j-2)}{j(j-1)(\mathrm
    N+j-1)(\mathrm N+j-2)} \\ & \stackrel{\mathrm N\to\infty}\sim 8
  \sum_{j=1}^\infty \frac{\mathrm N}{j(\mathrm
    N+j)^2} \stackrel{\mathrm
    N\to\infty}\sim 8 \sum_{j=1}^\infty \frac{1}{j(\mathrm N+j)} 
  \stackrel{\mathrm N\to\infty}\sim \frac{8}{\mathrm N}
  \sum_{j=1}^\infty \Big(\frac{1}{j} - \frac{1}{\mathrm N+j}\Big)
  \stackrel{\mathrm N\to\infty}{\sim} \frac{8\log\mathrm N}{\mathrm N}
\end{align*}
which proves the Proposition.
\end{proof}

\subsection{Numbers of ancestors near the tree top}
Let $u>0$. Define $S_u^{\mathrm N}$ to be the (random) number of
ancestors at time $-u$ in $\mathcal T^{\mathrm N}$ and $S_u$ the
number of ancestors by time $-u$ in $\mathcal T$, where $\mathcal
T^{\mathrm N}$ and $\mathcal T$ are defined as in Section
\ref{S:strong}. We give results on convergence of $S_u^{\mathrm N}$ as
$\mathrm N\to\infty$ (Lemma \ref{L:2}) and on convergence of $S_u$ as
$u\to 0$ (Lemma \ref{L:At}).
\label{sec:ancestors}

\begin{lemma}\label{L:2}
For $u>0$,
$$ S_u^{\mathrm N}\xrightarrow{\mathrm N\to\infty} S_u$$
almost surely and in $L^p$ for all $p>0$.
\end{lemma}

\begin{proof}
  With probability one, the completion of $R$ from \eqref{eq:R} is
  compact, so $\mathcal T$ \emph{comes down from infinity}, i.e.\ with
  probability one there are at most finitely many lines left by time
  $-u$. Since $S_u^{\mathrm N}$ is increasing with $\mathrm N$, the
  almost sure convergence follows. To see the $L^p$-convergence, note
  that all moments of the distribution of $S_u$ exist; see e.g.\
  \cite[Section 5.4]{Tavare1984}. Since $S_u^{\mathrm N}$ is bounded
  by $S_u$, the families $\big((S_u^{\mathrm N})^p\big)_{\mathrm
    N=1,2,...}$ are uniformly integrable and convergence in $L^p$
  follows.
\end{proof}

\begin{lemma}\label{L:At}
For the number of ancestors $S_u$
\begin{align}\label{eq:At2}
  u \cdot S_u \xrightarrow{u \to 0} 2,
\end{align}
almost surely and in $L^2$. Moreover,
\begin{align}
  \label{eq:At2b}
  \frac{S_u - 2/u}{\sqrt{ 2/{(3 u)}}} \xRightarrow{u\to 0} N(0,1).
\end{align}
In addition, for $u,v\to 0$, $u\leq v$, $u/v\to \Gamma \leq 1$
\begin{align}
  \label{eq:cltA2}
  \Big(\frac{S_u-2/u}{\sqrt{2/(3u)}},
  \frac{S_v-2/v}{\sqrt{2/(3v)}}\Big) \xRightarrow{\mbox{}} N(0,C)
\end{align}
with the covariance matrix $C$ given by
\begin{align}
  \label{eq:covMat}
  C = \left( \begin{matrix} 1 & \Gamma^{3/2} \\ \Gamma^{3/2} & 1 \end{matrix}
  \right).
\end{align}
\end{lemma}

\begin{remark}
  As a consequence of the previous lemma, the finite dimensional
  distributions of
  $\big((S_{tu}-\tfrac{2}{tu})/\sqrt{2/(3tu)}\big)_{t\geq 0}$ converge
  as $u\to 0$ to those of a Gaussian process $(A_t)_{t\geq 0}$ with
  covariance $\mathbb{COV}[A_s, A_t] = (s/t)^{3/2}$ for $s\leq t$.
\end{remark}

\begin{proof}[Proof of Lemma \ref{L:At}]
  The convergences \eqref{eq:At2} and \eqref{eq:At2b} can be found on
  p.\ 27 in \cite{Aldous1999}. For further use below, we prove
  \eqref{eq:At2b} in all detail. Define $T_n$ as the time it takes the
  coalescent to come down to $n$ lines, i.e.,
$$ T_n := \sum_{i=n+1}^\infty X_i.$$
Note that 
\begin{align}
  \label{eq:Kn1}
  \mathbb E[T_n] = \sum_{i=n+1}^\infty \frac{2}{i(i-1)} =
  \frac{2}{n},\qquad \mathbb V[T_n] = \sum_{i=n+1}^\infty
  \frac{4}{i^2(i-1)^2} \stackrel{n\to\infty}\sim \frac{4}{3n^3}.
\end{align}
The central ingredients in the proof are the two facts
\begin{equation}
  \label{eq:cltK}
  \begin{aligned}
    \mathbb P[S_u \leq n] & = \mathbb P[T_n \leq u],\\
    \frac{T_n - 2/n}{\sqrt{4/(3n^3)}} & \xRightarrow{n\to\infty}
    N(0,1),
  \end{aligned}
\end{equation}
where the second assertion is a consequence of \eqref{eq:Kn1} and
the central limit theorem. Hence we may define
\begin{align}
  \label{eq:asx}
  a_u(x) := \lfloor 2/u + x\sqrt{2/(3u )}\rfloor,
\end{align}
and write
\begin{align*}
  \mathbb P\Big[\frac{S_u - 2/u }{\sqrt{2/(3u )}}\leq x\Big] & =
  \mathbb P[S_u \leq a_u (x)] \\ & = \mathbb P[T_{a_u (x)} \leq u ]
  \\ & = \mathbb P\Big[\frac{T_{a_u (x)} - 2/a_u (x)} {\sqrt{4/(3a_u
      (x)^3)}} \leq \frac{u - 2/a_u (x)} {\sqrt{4/(3a_u
      (x)^3)}}\Big] \\ & \stackrel{u \to 0}{\sim} \mathbb
  P\Big[\frac{T_{a_u (x)} - 2/a_u (x)} {\sqrt{4/(3a_u (x)^3)}} \leq
  x \Big],
\end{align*}
since
\begin{align}
  \label{eq:aepsx}
  u  - 2/a_u (x) \stackrel{u \to 0}{\sim}
  x\sqrt{u ^3/6}, \qquad \sqrt{4/3(a_u (x)^3)}
  \stackrel{u \to 0}{\sim} \sqrt{u^3/6}.
\end{align}
Now, \eqref{eq:At2b} follows from \eqref{eq:cltK}. Since the event in
\eqref{eq:At2} is measurable with respect to the terminal
$\sigma$-algebra generated by the independent random variables $X_2,
X_3,...$, the convergence in \eqref{eq:At2} holds almost surely. The
$L^2$-convergence follows from moment results for $S_u$ given e.g.\
in \cite[Section 5.4]{Tavare1984}.

Let us turn to the proof of \eqref{eq:cltA2}. Since for $m\leq n$
$$ \mathbb{COV}[T_m,T_n] = \sum_{i=n+1}^\infty \mathbb V[X_i] 
\stackrel{n\to\infty}\sim \frac{4}{3n^3},$$ and using \eqref{eq:Kn1},
it is an easy exercise to show that for $m,n\to\infty, m/n\to
\Gamma\leq 1$,
\begin{align}
  \label{eq:TnTm}
  \Big(\frac{T_n-2/n}{\sqrt{4/(3n^3)}},
  \frac{T_m-2/m}{\sqrt{4/(3m^3)}}\Big) \Rightarrow N(0,C)
\end{align}
with the covariance matrix $C$ given in \eqref{eq:covMat}. To see
\eqref{eq:cltA2} from this, note that for $m\leq n$, $u\leq v$, an
extension of \eqref{eq:cltK} gives
\begin{align*}
  \mathbb P[T_n\leq u, T_m\leq v] & = \mathbb P[S_u\leq n, S_v\leq
  m]
\end{align*}
and thus for $x,y\in\mathbb R$, and using \eqref{eq:asx},
\begin{align*}
  \mathbb P\Big[ & \frac{S_u - 2/u}{\sqrt{2/(3u)}}\leq x,
  \frac{S_v-2/v}{\sqrt{2/(3v)}}\leq y\Big] = \mathbb P[S_u\leq
  a_{u}(x), S_v\leq a_{v}(y)] \\ & = \mathbb P[T_{a_u(x)}\leq u,
  T_{a_v(y)}\leq v] \\ & = \mathbb P\Big[ \frac{T_{a_u(x)}-
    2/a_u(x)}{\sqrt{4a_u(x)^3/3}} \leq
  \frac{u-2/a_u(x)}{\sqrt{4a_u(x)^3/3}}, \frac{T_{a_v(y)}-
    2/a_v(y)}{\sqrt{4a_v(y)^3/3}} \leq
  \frac{v-2/a_v(y)}{\sqrt{4a_v(y)^3/3}}\Big] \\ & \stackrel{u,v\to
    0}\sim \mathbb P\Big[ \frac{T_{a_u(x)}-
    2/a_u(x)}{\sqrt{4a_u(x)^3/3}} \leq x, \frac{T_{a_v(y)}-
    2/a_v(y)}{\sqrt{4a_v(y)^3/3}} \leq y\Big]
\end{align*}
by \eqref{eq:aepsx} and thus, if $u,v\to 0, u/v\to \Gamma\leq 1$,
\eqref{eq:cltA2} follows from \eqref{eq:TnTm} since $a_v(x)/a_u(y)
\stackrel{u,v\to 0}\sim u/v.$
\end{proof}

\subsection{The tree length near the tree top}
\label{S:44}
We analyze now the contribution to the tree length that comes from a
small time interval near the tree top. To this purpose we define
\begin{align}
\label{eq:treeLfi}
\Delta^{\mathrm N}_u := \int_0^u (S_v^{\mathrm N} - \mathbb
E[S_v^{\mathrm N}])dv, \qquad \Delta_u := \int_0^u (S_v - \mathbb
E[S_v])dv.
\end{align} Note that $\Delta_u^{\mathrm N}$ equals $\tilde
A_{0,u}^{\mathrm N}-\mathbb E[\tilde A_{0,u}^{\mathrm N}]$ from Figure
\ref{fig1a} in distribution. Again, we give results on convergence of
$\Delta_u^{\mathrm N}$ as $\mathrm N\to\infty$ (Lemma \ref{L:3}) and
of $\Delta_u$ as $u\to 0$ (Lemma \ref{L:treeLengthsSmallTimes}).

\begin{remark}[$\Delta_u$ as an $L^2$-limit]
  Since the integrand in the definition of $\Delta_u$ is unbounded,
  we have to make sure that the random variable $\Delta_u$
  exists. Indeed, using Lemma \ref{L:At} it is easy to check that
  $\big(\int_{2^{-n}}^u (S_v - \mathbb E[S_v])dv\big)_{n=1,2,...}$ is
  a Cauchy sequence in $L^2$, and we define $\Delta_u$ as its
  $L^2$-limit.  In particular, by continuity of ${\mathbb E}$ on $L^2$
  and Fubini's Theorem we further obtain ${\mathbb
    E}[\Delta_{u}]={\mathbb E}[\Delta_{u}^{\mathrm N}]=0$.
\end{remark}

\begin{lemma}
\label{L:3}
For the random variables $\Delta_u^{\mathrm N}$ and $\Delta_u$,
\begin{align}
  \label{eq:lengthSmallTimes1}
  \Delta_u^{\mathrm N}\xrightarrow{\mathrm N\to\infty}\Delta_u
\end{align}
in $L^2$.
\end{lemma}

\begin{proof}
  We start with proving the intuitively obvious fact that $S_v-
  S_v^{\mathrm N}$ and $ S_w- S_w^{\mathrm N}$ have nonnegative
  correlation. For $w\leq v$, we write
  \begin{equation}
    \label{eq:ass1}
    \begin{aligned}
      \mathbb{COV}[S_v & - S_v^{\mathrm N}, S_w- S_w^{\mathrm N}] =
      \mathbb{COV}\big[\mathbb E[S_v - S_v^{\mathrm N}|S_v,S_w],
      \mathbb E[S_w- S_w^{\mathrm N}|S_v,S_w]\big] \\ & \qquad \qquad
      \qquad \qquad \qquad \qquad \qquad \qquad + \mathbb E\big[
      \mathbb{COV}[S_v - S_v^{\mathrm N}, S_w- S_w^{\mathrm N}|S_v,
      S_w]\big] \\ & = \mathbb{COV}\Big[\frac{S_v(S_v-1)}{\mathrm
        N+S_v-1}, \frac{S_w(S_w-1)}{\mathrm N+S_w-1}\Big] + \mathbb
      E\Big[ \frac{S_v(S_v-1)\mathrm N(\mathrm N-1)}{(\mathrm N +
        S_w-1)(\mathrm N + S_v-1)(\mathrm N + S_v-2)}\Big],
    \end{aligned}
  \end{equation}
  where we have used \eqref{eq:distKi1} for the first and
  \eqref{eq:distKi2} for the second term. The second term on the
  r.h.s.\ is nonnegative, and so is the first term, since $i \mapsto
  \frac{i(i-1)}{\mathrm N+i-1}$ is increasing and $(S_v, S_w)$ are
  associated, i.e.
  \begin{align}\label{eq:ass}
    \mathbb{COV}[f(S_v), g(S_w)]\geq 0
  \end{align}
  for all non-decreasing functions $f, g$.  Indeed, to verify
  \eqref{eq:ass} it is enough to show this inequality for $f(S_v) =
  1_{S_v \ge \ell}$ and $g(S_w) = 1_{S_w \ge k}$. This, however, is
  clear since
  \begin{align*}
    \mathbb{COV}[1_{\{S_w\geq \ell\}}, 1_{\{S_v\geq k\}}] & =
    \mathbb{COV}[1_{\{S_w\geq \ell\}}, \mathbb E[1_{\{S_v\geq
      k\}}|S_w]] \geq 0
  \end{align*}
  by the well-known fact that a single random variable (here $S_w$) is
  associated and both $1_{\{S_w \ge \ell\}}$ and $\mathbb E[1_{\{S_v
    \ge k\}}|S_w]$ are non-decreasing functions of $S_w$. So we have
  proved that
  \begin{align}\label{eq:SvSw}
    \mathbb{COV}[S_v & - S_v^{\mathrm N}, S_w- S_w^{\mathrm N}] \geq
    0
  \end{align}
  for all $v,w\geq 0$.

  Now we come to the proof of \eqref{eq:lengthSmallTimes1}. By
  Fubini's Theorem and \eqref{eq:SvSw},
  \begin{align*}
    \mathbb E[(\Delta_u^{\mathrm N}-\Delta_u)^2] & = 2\int_0^u
    \int_0^v \mathbb E\big[\big(S_v^{\mathrm N}-S_v - \mathbb
    E[S_v^{\mathrm N}-S_v]\big) \big(S_w^{\mathrm N}-S_w - \mathbb
    E[S_w^{\mathrm N}-S_w]\big)\big]\,dw dv \\ & = 2\int_0^u \int_0^v
    \mathbb{COV}[S_v- S_v^{\mathrm N}, S_w- S_w^{\mathrm N}] dw dv \\
    & \leq 2\int_0^\infty \int_0^v \mathbb{COV}[S_v- S_v^{\mathrm N},
    S_w- S_w^{\mathrm N}] dw dv = \mathbb E[(\Delta_\infty^{\mathrm
      N}-\Delta_\infty)^2] \xrightarrow{\mathrm N\to\infty} 0
  \end{align*}
  by Proposition \ref{P:34}, since $(\Delta_\infty^{\mathrm
    N}-\Delta_\infty)^2$ is distributed as $B_{\mathrm N}^2$ in that
  Proposition.
\end{proof}

\begin{lemma}\label{L:treeLengthsSmallTimes}
For the random variables $\Delta_u$,
\begin{align}
  \label{eq:lengthSmallTimes2}
  \mathbb V[\Delta_u] \stackrel{u\to 0}{\sim} \tfrac 23 u.
\end{align}
\end{lemma}
\begin{proof}
  By Lemma \ref{L:At}, we see that for $w\leq v$
  \begin{align*}
    \mathbb{COV}[S_v,S_w] \stackrel{w\leq v\to 0}\sim
    \sqrt{\frac{2}{3w}\frac{2}{3v}\frac{w^3}{v^3}} = \frac 23
    \frac{w}{v^2}.
  \end{align*}
  Hence, by Fubini's Theorem,
  $$ \mathbb V[\Delta_u] =
  2\int_0^u\int_0^v \mathbb{COV}[S_w,S_v] dw dv \stackrel{u\to 0}\sim
  \frac 43 \int_0^u \int_0^v \frac{w}{v^2} dw dv = \frac 23 u.$$
\end{proof}
\begin{remark}
Note that $\Delta_u$ is an integral over approximately Gaussian
random variables for small $u$. In addition, for $s\leq t$,
\begin{align*}
  \mathbb{COV}\Big[ & \frac{\Delta_{su}}{\sqrt{\tfrac 23 su}},
  \frac{\Delta_{tu}}{\sqrt{\tfrac 23 tu}}\Big] =
  \frac{3}{2u\sqrt{st}} \int_0^{su} \int_0^{tu} \mathbb{COV}[S_v, S_w]
  dw dv \\ & \stackrel{u\to 0}\sim \frac{3}{2u\sqrt{st}} \Big( 2
  \int_0^{su} \int_0^v \frac 23 \frac{w}{v^2} dw dv + \frac 23
  \int_0^{su}wdw \int_{su}^{tu} \frac{1}{v^2} dv\Big) \\ & =
  \sqrt{\frac{s}{t}} + \frac{1}{2}\Big(\sqrt{\frac{s}{t}} -
  \sqrt{\frac{s^3}{t^3}}\Big) = \frac 32 \sqrt{\frac s t} - \frac 12
  \sqrt{\frac{s^3}{t^3}}.
\end{align*}
Hence, as an extension of \eqref{eq:lengthSmallTimes2}, we see that
the finite dimensional distributions of $\Big(
{\Delta_{tu}}/{\sqrt{\tfrac 23 tu}}\Big)_{t\geq 0}$ converge as $u\to
0$ to those of a centered Gaussian process $(A_t)_{t\geq 0}$ with
covariance $\mathbb{COV}[A_s, A_t] = \frac 32 \sqrt{\frac s t} - \frac
12 \sqrt{\frac{s^3}{t^3}}$ for $s\leq t$.
\end{remark}

\subsection{The evolution of the $f$ oldest families in the Moran
  model}
\label{Moran}
Consider the graphical representation of a Moran model given in Figure
\ref{fig1}. For any time $t\in\mathbb R$ the tree $\mathscr T^{\mathrm
  N}_t$ can be identified with a random subset of $(-\infty;t]\times
\{1,...,\mathrm N\}$ which we continue to denote by $\mathscr
T^{\mathrm N}_t$.  For any $t\in\mathbb R$ we define the
inter-coalescence times
\begin{align}\label{eq:Mor1}
  X_i^{\mathrm N}(t) := \text{ length of the time interval in which
    $\mathscr T_t^{\mathrm N}$ has $i$ lines}
\end{align}
for $i=2,...,\mathrm N$. We denote by $\mathscr R^{\mathrm N}$ the
rate $\binom{\mathrm N}2$ Poisson process of all resampling events,
viewed as a random subset of $\mathbb R$. For $f=2,...,\mathrm N$ and
$t\in\mathbb R$, consider the $f$ oldest families at time $t$,
i.e. the $f$ subtrees of $\mathscr T^{\rm N}_t$ (call them $F^{\rm
  N}_t(1), \ldots, F^{\rm N}_t(f)$) whose union is $\mathscr T_t^{\rm
  N} \cap\big(\{[t-\sum_{i=f+1}^{\rm N}X_j^{\rm N};
t]\times\{1,\ldots, \rm N\}\big)$. Let
$$\mathscr D^{{\mathrm N}, f}:=\{t: X_f^{\mathrm N}(t) \neq X_f^{\mathrm N}(t-)\}
\subseteq \mathscr R^{\mathrm N},$$ i.e.\ $\mathscr D^{{\rm N},f}$ is
the point process of times when one of the $f$ oldest families gets
extinct. We set $Z^{\mathrm N, f}_t = (Z^{\mathrm N, f}_t(1), ...,
Z^{\mathrm N, f}_t(f))$, where
$$Z_t^{{\rm N},f}(i):= \mbox{ size of the family $F^{\rm N}_t(i)$},$$
i.e.\ the number of leaves in the tree $F^{\rm N}_t(i)$.  Note that
$Z^{{\mathrm N}, f}_t \in \mathbb Z^f_{++\mathrm N} :=
\{(z_1,\ldots,z_f): z_i \ge 1, z_1+\ldots+z_f=\mathrm N\}$. The
following lemma is essential in the proof of Theorem \ref{TH1}.

\begin{lemma}[The $f$ oldest families in a Moran model]
  \label{L:moran}
  Fix $f\in\{2,...,\mathrm{N}\}$. Let $\mathscr R^{\mathrm N} =
  \{\tau_n: n\in\mathbb Z\}$ with $...<\tau_0<\tau_1<...$.
  \begin{enumerate}
  \item The uniform distribution on $\mathbb Z^f_{++\mathrm N}$ is the
    stationary distribution for the Markov chains $(Z^{\mathrm
      N,f}_t)_{t\in\mathbb R}$ and $(Z^{\mathrm
      N,f}_{\tau_n})_{n\in\mathbb Z}$.
  \item The events $(\{\tau_n\in\mathscr D^{\mathrm
      N,f}\})_{n\in\mathbb Z}$ are independent and $\mathbb
    P[\tau_n\in\mathscr D^{\mathrm N,f}] = \binom f 2/\binom{\mathrm
      N}2$. In particular, $\mathscr D^{\mathrm N,f}$ is a Poisson
    process with rate $\binom f2$.
  \end{enumerate}
\end{lemma}
\begin{proof}
  Let $t\in\mathbb R$. Looking at the Moran genealogy during the time
  interval $\big[t-\sum_{i=f+1}^{\mathrm N} X_i^{\mathrm N}(t);
  t\big]$, one sees that the $f$ oldest families have been built up
  through a P\'olya urn with $f$ ancestral balls. This already
  explains 1.\ for a fixed (non-random) $t$. The same arguments apply
  if $t$ is a resampling time. Hence, we have proved assertion 1.

  Let us now consider the dynamics of the sizes of the $f$ oldest
  families $(Z^{\mathrm N,f}_{\tau_n})_{n\in\mathbb Z}$ in its
  equilibrium, the uniform distribution on $\mathbb Z^f_{++\mathrm
    N}$. The Moran model is set up as follows: Each pair of balls is
  chosen at unit rate. As soon as a pair is chosen, one of the two
  balls (selected at random from the pair) is transferred into the box
  of the other ball. (If two balls from the same box are chosen,
  nothing changes.)
  For the transition of the Markov chain, we compute for
  $(z_1,...,z_f)\in\mathbb Z^f_{++\mathrm N} $
  \begin{align*}
    \mathbb P & [Z^{\mathrm N,f}_{\tau_{n+1}} =
    (z_1,...,z_f),\tau_{n+1}\notin \mathscr D^{\mathrm N,f}] =
    \sum_{{i,j=1}\atop{i\neq j}}^f \mathbb P[Z^{\mathrm
      N,f}_{\tau_{n}} = (z_1,...,z_i+1,...,z_j-1,...,z_f)]
    \frac{\tfrac 12(z_i+1)(z_j-1)}{\binom{\mathrm N}{2}} \\ & \qquad
    \qquad \qquad \qquad \qquad \qquad \qquad \qquad + \sum_{i=1}^f
    \mathbb P[Z^{\mathrm N,f}_{\tau_{n}} =
    (z_1,...,z_i,...,z_j,...,z_f)]
    \frac{\binom{z_i}{2}}{\binom{\mathrm N}{2}} \\ & = \mathbb P
    [Z^{\mathrm N,f}_{\tau_{n}} = (z_1,...,z_f)] \cdot
    \frac{1}{\binom{\mathrm N}{2}} \cdot \tfrac 12 \cdot \Big(\Big(
    \sum_{i=1}^f z_i\Big)^2 + \sum_{{i,j=1}\atop{i\neq j}}^f
    (z_i-z_j-1) - \sum_{i=1}^f z_i\Big) \\ & = \Big( 1 -
    \frac{\binom{f}{2}}{\binom{\mathrm N}{2}}\Big) \cdot\mathbb P
    [Z^{\mathrm N,f}_{\tau_{n}} = (z_1,...,z_f)].
  \end{align*}
  This calculation reveals two things: first, by summing over all
  possible $(z_1,...,z_f)$ on both sides, we see that $\mathbb
  P[\tau_{n+1}\in \mathscr D] = \binom{f}{2}/\binom{\mathrm
    N}{2}$. Secondly, we see that $Z^{\mathrm N,f}_{\tau_{n+1}}$, given
  $\tau_{n+1}\in \mathscr D$, is again uniformly distributed on
  $\mathbb Z^f_{++\mathrm N}$. The last assertion implies that
  $(\{\tau_n\in\mathscr D^{\mathrm N,f}\})_{n\in\mathbb Z}$ is
  independent, which finishes the proof of assertion 2.


\end{proof}

\section{Proof of Theorem \ref{TH1}}
\label{S:Sec4}
We will prove Theorem \ref{TH1} in three steps. For convergence of
$\mathscr L^{\mathrm N}$ we need to show (see e.g. \cite[Lemma
3.4.3]{EthierKurtz1986})
\begin{enumerate}
\item[(a)] The sequence of processes $\mathscr L^{\mathrm N}$ is tight
 in $\mathbb D$,
\item[(b)] The finite-dimensional distributions of the sequence
$\mathscr L^{\mathrm N}$ converge.
\end{enumerate}

The main work is to show tightness of $(\mathscr L^{\mathrm
  N})_{\mathrm N\in\mathbb N}$ in $\mathbb D$. For this, it is enough
to show (see \cite[Theorem 3.8.6 and Theorem 3.8.8]{EthierKurtz1986})
\begin{align}
\label{eq:oneDim} 
(\mathscr L^{\mathrm N}_0)_{\mathrm N\in\mathbb N} \text{ is tight
  in } \mathbb R
\end{align}
and there exists $\beta>0$ and $\theta>1$ such that\footnote{For
 functions $a, b: \mathbb R_+ \to \mathbb R_+$ we write $a(h)
 \lesssim b(h)$ iff there is a $C>0$, independent of any other
 parameter, such that $a(h)\leq C \cdot b(h)$ for all $h>0$.} for all
$t\in\mathbb R$,
\begin{align}
\label{eq:tight2}
\limsup_{\mathrm N\to\infty} \mathbb E[1\wedge |\mathscr L^{\mathrm
  N}_{t+h}-\mathscr L^{\mathrm N}_t|^\beta \wedge |\mathscr L^{\mathrm
  N}_{t}-\mathscr L^{\mathrm N}_{t-h}|^\beta ] \lesssim h^\theta.
\end{align}
Proposition \ref{L:Gumbel} already shows that $\tfrac 12 \mathscr
L_0^{\mathrm N}$ converges to a Gumbel distributed random variable,
which implies \eqref{eq:oneDim}. Now, the main work is to show
\eqref{eq:tight2}, which will be done in Step 1. In Step 2 we show
convergence of finite-dimensional distributions. Step 3 then shows
\eqref{Thm1}.

\begin{step}[Proof of \eqref{eq:tight2}]
  Consider a Moran model of size $\mathrm N$.  We will use the same
  notation as in Subsection \ref{Moran}. For any two time points $t$
  and $t+h$, we put
  $$ F^{\mathrm N}_{(t,t+h]} := \min\{i=2,...,\mathrm N: 
  X_i^{\mathrm N}(t) \neq X_i^{\mathrm N}(t+h)\}.$$ Note that
  $\{F^{\mathrm N}_{(t,t+h]} \geq f\}$ is the event that none of the
  $f$ oldest families gets extinct during $(t; t+h]$. Since the point
  process of losses of one of the $f$ oldest families is Poisson with
  rate $\binom f2$ by Lemma \ref{L:moran}, we conclude that
  \begin{align}\label{eq:ffour}
    \mathbb P[F^{\mathrm N}_{(t-h,t]} \vee F^{\mathrm N}_{(t,t+h]} <
    f] & = (1-e^{-\binom f 2 h})^2 \lesssim f^4 h^2 \wedge 1.
  \end{align}
  Using this equation, we will now show \eqref{eq:tight2} for $\theta
  = \tfrac{10}9$ and $\beta=10$. We write
   \begin{equation}
     \label{eq:p11}
     \begin{aligned}
       \mathbb E & [1\wedge (\mathscr L^{\mathrm N}_{t+h}-\mathscr
       L^{\mathrm N}_t)^{10} \wedge (\mathscr L^{\mathrm
         N}_{t}-\mathscr L^{\mathrm N}_{t-h})^{10} ] \\ & \leq
       \mathbb E [1\wedge (\mathscr L^{\mathrm N}_{t+h}-\mathscr
       L^{\mathrm N}_t)^{10}; F^{\mathrm N}_{(t,t+h]} \geq F^{\mathrm
         N}_{(t-h,t]}] + \mathbb E [1\wedge (\mathscr L^{\mathrm
         N}_{t}-\mathscr L^{\mathrm N}_{t-h})^{10}; F^{\mathrm
         N}_{(t-h,t]} \geq F^{\mathrm N}_{(t,t+h]}].
     \end{aligned}
   \end{equation}
   We will next bound the first term on the right hand; the bound for
   the second term is obtained in the same manner. We get
   \begin{equation}
     \label{eq:p12}
     \begin{aligned}
       \mathbb E & [1\wedge (\mathscr L^{\mathrm N}_{t+h}-\mathscr
       L^{\mathrm N}_t)^{10}; F^{\mathrm N}_{(t,t+h]} \geq F^{\mathrm
         N}_{(t-h,t]}] \\ & \leq \mathbb P[F^{\mathrm N}_{(t-h,t]}
       \leq F^{\mathrm N}_{(t,t+h]} \leq h^{-2/9} ] + \mathbb E [
       (\mathscr L^{\mathrm N}_{t+h}-\mathscr L^{\mathrm N}_t)^{10};
       F^{\mathrm N}_{(t,t+h]} \geq h^{-2/9}] \\ & \leq \mathbb
       P[F^{\mathrm N}_{(t-h,t]} \vee F^{\mathrm N}_{(t,t+h]} \leq
       h^{-2/9} ] \\ & \qquad \qquad \qquad + \mathbb E\Big[
       \Big( \sum_{i=2}^{\mathrm N} i\Big( X_i^{\mathrm N}(t+h) -
       \tfrac{1}{\binom{i}{2}}\Big) - \sum_{i=2}^{\mathrm N} i\Big(
       X_i^{\mathrm N}(t) - \tfrac{1}{\binom{i}{2}}\Big)\Big)^{10};
       F^{\mathrm N}_{(t,t+h]} \geq h^{-2/9}\Big] \\ & \lesssim
       h^{10/9} + \mathbb E\Big[ \Big( \sum_{i=\lfloor
         h^{-2/9}\rfloor}^{\mathrm N} i\Big( X_i^{\mathrm N}(t+h) -
       \tfrac{1}{\binom{i}{2}}\Big) - \sum_{i=\lfloor
         h^{-2/9}\rfloor}^{\mathrm N} i\Big( X_i^{\mathrm N}(t) -
       \tfrac{1}{\binom{i}{2}}\Big)\Big)^{10}\Big]\, ,
     \end{aligned}
   \end{equation}
   where the last inequality follows from \eqref{eq:ffour} and the
   fact that $X_i(t+h)=X_i(t)$ on $\{ F^{\mathrm N}_{(t,t+h]} \geq
   h^{-2/9}\}$ for all $i<\lfloor h^{-2/9}\rfloor$ by definition of
   $F^{\mathrm N}_{(t,t+h]}$.  It remains to bound the second term in
   the last line. For this, we define for $n=1,2,...$
   \begin{align*}
     a_n(h) := \sum_{i=\lfloor h^{-2/9}\rfloor}^{\mathrm N} \mathbb E
     \Big[ \Big( i\Big( X_i^{\mathrm N}(t) -
     \tfrac{1}{\binom{i}{2}}\Big) \Big)^{n}\Big] 
   \end{align*}
   and observe that
   \begin{equation}
     \label{eq:p13}
     \begin{aligned}
       a_1(h) = 0, \qquad \qquad a_n(h) \lesssim h^{2(n-1)/9},
     \end{aligned}
   \end{equation}
   since the $n$th central moment of an exponentially distributed
   random variable with parameter $\lambda$ is proportional to
   $\lambda^{-n}$. In addition, we use that $(x-y)^n \leq (2x)^n +
   (2y)^n$ for even $n$ and all $x,y\in\mathbb R$, and independence
   of $X_2^{\mathrm N}(t),...,X_{\mathrm N}^{\mathrm N}(t)$ as well as of
   $X_2^{\mathrm N}(t+h),...,X_{\mathrm N}^{\mathrm N}(t+h)$ to obtain
   \begin{equation}
     \label{eq:p14}
     \begin{aligned}
       \mathbb E \Big[ \Big( \sum_{i=\lfloor
         h^{-2/9}\rfloor}^{\mathrm N} & i\Big( X_i^{\mathrm N}(t+h) -
       \tfrac{1}{\binom{i}{2}}\Big) - \sum_{i=\lfloor
         h^{-2/9}\rfloor}^{\mathrm N} i\Big( X_i^{\mathrm N}(t) -
       \tfrac{1}{\binom{i}{2}}\Big)\Big)^{10}\Big] \\ & \lesssim
       \mathbb E \Big[ \Big( \sum_{i=\lfloor
         h^{-2/9}\rfloor}^{\mathrm N} i\Big( X_i^{\mathrm N}(t+h) -
       \tfrac{1}{\binom{i}{2}}\Big) \Big)^{10}\Big] + \mathbb E \Big[
       \Big( \sum_{i=\lfloor h^{-2/9}\rfloor}^{\mathrm N} i\Big(
       X_i^{\mathrm N}(t) - \tfrac{1}{\binom{i}{2}}\Big)
       \Big)^{10}\Big] \\ & \lesssim \sum_{k=1}^{10}
       \sum_{{n_1,...,n_k} \atop {n_1+...+n_k=10}}^{10}
       a_{n_1}(h)\cdots a_{n_k}(h) \lesssim (a_2(h))^5 \lesssim
       h^{10/9}.
     \end{aligned}
   \end{equation}
   by \eqref{eq:p13}. Plugging \eqref{eq:p14} and \eqref{eq:p12} into
   \eqref{eq:p11} shows that \eqref{eq:tight2} holds with $\beta=10$
   and $\theta = \tfrac{10}{9}$.
 \end{step}

\begin{step}[Convergence of finite-dimensional distributions]
 Fix $t_1<...<t_n$. We will show that $(\mathscr L_{t_1}^{\mathrm
   N},..., \mathscr L_{t_n}^{\mathrm N})$ converges weakly for
 $\mathrm N\to\infty$. The strategy is to define a probability space
 on which all $\mathscr L_{t_i}^{\mathrm N}$,  $i=1,...,n,
 \mathrm N=2,3,...$, are defined.

 Consider a coalescent, started with infinitely many lines which are
 numbered by $1_n,2_n,...$ for some time $t_n-t_{n-1}$. Denote the
 number of the ancestors at time $t_n-t_{n-1}$ by $S_n$ . Number the
 lines going back from these by $1_{n-1},..., (S_n)_{n-1}$ and augment
 them by lines numbered $(S_n+1)_{n-1}, (S_n+2)_{n-1},...$. Let these
 infinitely many lines coalesce for some time $t_{n-1}-t_{n-2}$,
 number the $S_{n-1}$ ancestors at time $t_{n-1}-t_{n-2}$ by
 $1_{n-2},..., (S_{n-1})_{n-2}$ and augment their lines by lines
 numbered $(S_{n-1}+1)_{n-2}, (S_{n-1}+2)_{n-2},...$. In this way we
 get iteratively $n$ genealogies for an infinite population back from
 times $t_1,...,t_n$. Considering the compensated tree lengths of
 lines numbered $1_i,...,\mathrm N_i$ gives the compensated tree
 length of a population of size $\mathrm N$ at time $t_i$,
 $i=1,...,n$. Moreover, as shown in Proposition \ref{L:Gumbel2}, these
 tree lengths converge in $L^2$ as $\mathrm N\to\infty$ for each
 $i=1,...,n$. Since $L^2$-convergence implies convergence in
 probability, which, in turn, implies weak convergence, we are done.
\end{step}

\begin{step}[Decomposition of $\mathscr L_t^{\mathrm N} - \mathscr
 L_0^{\mathrm N}$ and proof of \eqref{Thm1}]
 Recall the graphical representation of a Moran model from Figure
 \ref{fig1}. Using the random set $\mathscr T^{\mathrm N}_t \subseteq
 (-\infty,t]\times \{1,...,\mathrm N\}$ (recall Subsection
 \ref{Moran}), we have the representation $(\mathscr L^{\mathrm
   N}_t)_{t\in\mathbb R} \stackrel{d}{=} (\lambda^{\mathrm N}(\mathscr
 T^{\mathrm N}_t)-2 \log \mathrm N)_{t\in\mathbb R}$ where
 $\lambda^{\mathrm N}$ is Lebesgue measure on $\mathbb R\times
 \{1,...,\mathrm N\}$. We set
\begin{equation}
  \label{eq:step21}
  \begin{aligned}
    \tilde A^{\mathrm N}_{0,t} & := \lambda^{\mathrm N} (\mathscr
    T^{\mathrm N}_t \setminus \mathscr T^{\mathrm N}_0), \qquad
    A^{\mathrm N}_{0,t} := \tilde A^{\mathrm N}_{0,t} - \mathbb
    E[\tilde
    A^{\mathrm N}_{0,t}],\\
    \tilde B^{\mathrm N}_{0,t} & := \lambda^{\mathrm N} (\mathscr
    T^{\mathrm N}_0 \setminus \mathscr T^{\mathrm N}_t), \qquad
    B^{\mathrm N}_{0,t} := \tilde B^{\mathrm N}_{0,t} - \mathbb
    E[\tilde B^{\mathrm N}_{0,t}],
  \end{aligned}
\end{equation}
compare also with Figure \ref{fig1a}.  Note that $\mathbb
E[\lambda^{\mathrm N} (\mathscr T^{\mathrm N}_t \setminus \mathscr
T^{\mathrm N}_0)]=\mathbb E[\lambda^{\mathrm N} (\mathscr T^{\mathrm
  N}_0 \setminus \mathscr T^{\mathrm N}_t)]$ due to stationarity, and
thus
\begin{align}\label{eq:Lrepr}
  \mathscr L_t^{\mathrm N} - \mathscr L_0^{\mathrm N}
  \stackrel{d}{=} A^{\mathrm N}_{0,t} - B^{\mathrm N}_{0,t}.
\end{align}
For the infinitesimal variance, we find by the convergence of finite
dimensional distributions and \eqref{eq:Lrepr} that
\begin{align}\label{ANBN} \mathbb E[(\mathscr L_t - \mathscr L_0)^2] = \lim_{\mathrm N\to\infty} 
  \mathbb E[(\mathscr L_t^{\mathrm N} - \mathscr L_0^{\mathrm N})^2] =
  \lim_{\mathrm N\to\infty} \mathbb E[(A_{0,t}^{\mathrm N} -
  B_{0,t}^{\mathrm N})^2].
\end{align}
 From \eqref{eq:step21} and
\eqref{eq:treeLfi} we conclude 
that $A_{0,t}^{\mathrm N} \stackrel d = \Delta_t^{\mathrm N}$, again see
 Figure \ref{fig1a}. We have, using the $L^2$-convergence from
Lemma \ref{L:3} and Lemma \ref{L:treeLengthsSmallTimes}
\begin{align} \label{AN} \lim_{\mathrm N\to\infty} \mathbb V[A_{0,t}^{\mathrm N}] = \lim_{\mathrm
  N\to\infty} \mathbb V[ \Delta_t^{\mathrm N}] = \mathbb
V[\Delta_t] \stackrel{t\to 0}\sim \tfrac 23
t.
\end{align}
For the variance of $B_{0,t}^{\mathrm N}$, note that $B_{0,t}^{\mathrm
  N} \stackrel{d}=B_{\mathrm N} - B_{S_t^{\mathrm N}}$, where
$B_{\mathrm N}$ and $S_t^{\mathrm N}$ are as in Sections
\ref{sec:subtrees} and \ref{sec:ancestors}, with $\mathcal T :=
\mathscr T_0$ and $S_t^{\mathrm N}$ independent. We thus have for
fixed $t$ and $N \to \infty$, because of Proposition \ref {P:34},
\begin{align}\label{BN}
  \mathbb V[B_{0,t}^{\mathrm N}] = \mathbb V[B_{\mathrm N} -
  B_{S_t^{\mathrm N}}] \stackrel{\mathrm N\to\infty}\sim \mathbb
  V[B_{S_t^{\mathrm N}}] \stackrel{\mathrm N\to\infty}\sim \mathbb
  V[B_{S_t}].
\end{align}
Since $t\cdot S_t \to 2$ almost surely as $t\to 0$ (see Lemma
\ref{L:At}), we conclude from \eqref{BN} and Proposition \ref{P:34},
\eqref{eq:toshow1}, that
\begin{align}\label{BN1}
  \lim_{\mathrm N\to\infty} \mathbb V[B_{0,t}^{\mathrm N}] = \mathbb
  V[B_{S_t}] \stackrel{t\to 0}\sim \mathbb V[B_{\lfloor \tfrac
    2t\rfloor}] \stackrel{t\to 0}\sim 4t|\log t|\,.
\end{align}
Finally, combining \eqref{ANBN}, \eqref{AN} and \eqref {BN1}, and
noting that $\mathbb{COV}[A_{0,t}^{\mathrm N}, B_{0,t}^{\mathrm N}]
\stackrel{t\to 0}\lesssim t\sqrt{\tfrac 83|\log t|}$ by the
Cauchy-Schwartz inequality, we arrive at
\begin{align*}
  \mathbb E[(\mathscr L_t - \mathscr L_0)^2] 
  \stackrel{t \to 0}\sim 4t \, |\log t|.
\end{align*}
This completes the proof of Theorem \ref{TH1}.
\end{step}

\section{Proofs of strong convergence results}
\label{S:proofs2} In this section we prove Propositions \ref{L:Gumbel2} and
 \ref{P:TH1}.

\subsection{Proof of Proposition \ref{L:Gumbel2}}
From Proposition \ref{L:Gumbel} we know that $\Lambda^{\mathrm N}_1$
converges weakly as $\mathrm N\to\infty$ to a random variable
$\Lambda$ such that $\tfrac 12 \Lambda$ is Gumbel distributed. Since
$\Lambda^{\mathrm N}_1$ is a sum of independent random variables,
Kolmogorov's three series criterion shows that the convergence holds
almost surely as well. Moreover, since seconds moments converge in
\eqref{eq:P11}, the convergence also holds in $L^2$.

Next, we will show that
$$ \Lambda^{\mathrm N}_1 - \Lambda^{\mathrm N}_2 \xrightarrow{\mathrm N\to\infty} 0$$
in $L^2$. Together with the $L^2$-convergence of $\Lambda^{\mathrm
  N}_1$ this gives \eqref{eq:P12}. We compute directly, recalling
\eqref{eq:lambda2} and using that $\Lambda^{\mathrm N}_1 \stackrel d =
\Lambda^{\mathrm N}_2$ in the third equality and Lemma \ref{L1},
\begin{align*}
  \mathbb E\big[\big( \Lambda^{\mathrm N}_1 - \Lambda^{\mathrm
    N}_2\big)^2\big] & = \mathbb E\Big[\Big( \sum_{i=2}^\infty
  K_i^{\mathrm N} X_i - \sum_{i=2}^{\mathrm N} i X_i\Big)^2\Big] \\
  & = \mathbb E\Big[\Big( \sum_{i=2}^\infty K_i^{\mathrm N}
  X_i\Big)^2\Big] + \mathbb E\Big[\Big( \sum_{i=2}^{\mathrm N} i
  X_i\Big)^2\Big] - 2 \mathbb E\Big[\sum_{i=2}^\infty
  \sum_{j=2}^{\mathrm N} K_i^{\mathrm N} j X_i X_j\Big] \\ & = 2
  \Big(\mathbb E\Big[\Big( \sum_{i=2}^{\mathrm N} i X_i\Big)^2\Big] -
  \sum_{i=2}^\infty \sum_{j=2}^{\mathrm N} \frac{ij\mathrm N}{\mathrm
    N+i-1} \mathbb E[X_i X_j]\Big) \\ & = 2\Big( \mathbb
  V\Big[\sum_{i=2}^{\mathrm N}iX_i\Big] + \Big(\mathbb E\Big[
  \sum_{i=2}^{\mathrm N}iX_i\Big]\Big)^2 - \sum_{i=2}^\infty
  \sum_{j=2}^{\mathrm N}\frac{ij\mathrm N}{\mathrm
    N+i-1}\frac{4}{i(i-1)j(j-1)}(1+\delta_{ij})\Big) \\ & = 2\Big( 4
  \sum_{i=1}^{\mathrm N-1}\frac{1}{i^2} + \Big(2\sum_{i=1}^{\mathrm
    N-1}\frac 1i \Big)^2 - 4 \sum_{j=1}^{\mathrm N-1}\frac{1}{j}
  \sum_{i=1}^\infty \frac{\mathrm N}{i(\mathrm N+i)} - 4
  \sum_{i=1}^{\mathrm N-1} \frac{\mathrm N}{i^2(\mathrm N+i)} \Big)\\
  & = 8 \sum_{i=1}^{\mathrm N-1}\frac{1}{i^2}\Big( 1 - \frac{\mathrm
    N}{\mathrm N+i}\Big) + 8 \Big(\sum_{i=1}^{\mathrm N-1}\frac 1i
  \Big)^2 - 8 \sum_{j=1}^{\mathrm N-1}\frac 1j \sum_{i=1}^{\infty}
  \Big(\frac{1}{i} - \frac{1}{\mathrm N+i}\Big) \\ & = 8
  \sum_{i=1}^{\mathrm N-1}\frac{1}{i(\mathrm N+i)} = \frac{8}{\mathrm
    N}\sum_{i=1}^{\mathrm N-1} \Big(\frac 1i - \frac{1}{\mathrm
    N+i}\Big) \stackrel{\mathrm N\to\infty} \sim \frac{8 \log\mathrm
    N}{\mathrm N} \xrightarrow{\mathrm N\to\infty} 0.
\end{align*}

\subsection{Proof of Proposition \ref{P:TH1}}
We recall \cite[Lemma A2.1]{DonnellyKurtz1996}:

\begin{proposition}\label{P:tom}
  Let $(X^n)_{n=1,2,...}$ be a sequence of processes with sample paths
  in $\mathbb D$, defined on the same probability space. Suppose that
  $(X^n)_{n=1,2,...}$ is relatively compact in $\mathbb D$ (in the
  sense of convergence in distribution) and that for a dense set
  $H\subseteq\mathbb R$, $(X^n_t)_{n=1,2,...}$ converges in
  probability in $\mathbb R$ for each $t\in H$. Then, there is a
  process $X$ such that $d_{\mathrm
    {Sk}}(X^n,X)\xrightarrow{n\to\infty}0$ in probability.
\end{proposition}

We use this Proposition for $(\mathscr L^{\mathrm {ld},\mathrm N})_{\mathrm
 N=2,3,...}$. First, $\mathscr L^{\mathrm {ld},\mathrm N} \stackrel d =
\mathscr L^{\mathrm N}$ with $\mathscr L^{\mathrm N}$ as in Theorem
\ref{TH1}. Hence, as Theorem \ref{TH1} shows, $(\mathscr L^{\mathrm {ld},\mathrm
 N})_{\mathrm N=2,3,...}$ converges weakly. In particular, the
sequence is relatively compact in $\mathbb D$.

For all $t\in\mathbb R$, we have that $\mathscr T^{\mathrm {ld},\mathrm N}
\stackrel d = \mathscr T^{\mathrm N}$ with $\mathscr T^{\mathrm N}$
from \eqref{eq:defNat}. Consequently, $\mathscr L^{\mathrm {ld},\mathrm N}
\stackrel d = \Lambda_2^{\mathrm N}$ and there exists a random
variable $\mathscr L^{\mathrm {ld}}_t$ such that $\mathscr L^{\mathrm {ld}, \mathrm N}_t -
\mathscr L^{\mathrm {ld}}_t \xrightarrow{\mathrm N\to\infty} 0 $ in $L^2$. Since
the $L^2$-convergence implies convergence in probability we have
proved Proposition \ref{P:TH1}.

\subsubsection*{Acknowledgments}
We thank Tom Kurtz for teaching us Proposition \ref{P:tom}, and Martin
M\"ohle for reminding us of Kingman's notions of temporal and natural
coupling. Also, we thank an anonymous referee for helpful remarks.

PP and AW obtained travel support from the DFG, Bilateral Research
Group FOR 498. PP and HW are supported by the BMBF, Germany, through
FRISYS (Freiburg Initiative for Systems biology), Kennzeichen 0313921.

\bibliographystyle{plain}


\begin{thebibliography}{10}

\bibitem{Aldous1999}
D.~Aldous.
\newblock Deterministic and stochastic models for coalescence (aggregation and
  coagulation): a review of the mean-field theory for probabilists.
\newblock {\em Bernoulli} 5(1):3--48, 1999.

\bibitem{BBL2010}
J.~Berestycki, N.~Berestycki, and V.~Limic.
\newblock {The $\Lambda$-coalescent speed of coming down from infinity }.
\newblock {\em Ann. Probab.} 38(1):207--233, 2010.

\bibitem{BBS2009}
J.~Berestycki, N.~Berestycki, and J.~Schweinsberg.
\newblock Small-time behavior of beta coalescents.
\newblock {\em Ann. Inst. H. Poincaré Probab. Statist.} 44(2):214--238, 2008.

\bibitem{Berest09}
N.~Berestycki.
\newblock Recent progress in coalescent theory.
\newblock {\em Ensaios Matematicos} 16:1--193, 2009.

\bibitem{BlumFrancois2005}
M.~G. Blum and O.~Francois.
\newblock Minimal clade size and external branch length under the neutral
  coalescent.
\newblock {\em Adv. Appl. Probab.} 37:647--662, 2005.

\bibitem{CaliebeNeiningerKrawczakRoesler2007}
A.~Caliebe, R.~Neininger, M.~Krawczak, and U.~R{\"o}sler.
\newblock On the length distribution of external branches in coalescence trees:
  genetic diversity within species.
\newblock {\em Theo. Popul. Biol.} 72(2):245--252, 2007.

\bibitem{Dawson1993}
D.A. Dawson.
\newblock Measure-valued {M}arkov processes.
\newblock In P.L. Hennequin, editor, {\em \'Ecole d'\'Et\'e de Probabilit\'es
  de Saint-Flour XXI--1991}, volume 1541 of {\em Lecture Notes in Mathematics},
  pages 1--260, Berlin, 1993. Springer.

\bibitem{DelmasEtAl2008}
J.-F. Delmas, J.-S. Dhersin, and A.~Siri-Jegousse.
\newblock Asymptotic results on the length of coalescent trees.
\newblock {\em Ann. Appl. Prob.} 18(3):997--1025, 2008.

\bibitem{DelmasEtAl2010}
J.-F. Delmas, J.-S. Dhersin, and A.~Siri-Jegousse.
\newblock {On the two oldest families for the Wright-Fisher process}.
\newblock {\em Electron. J. Probab.} to appear, 2010.

\bibitem{DonnellyKurtz1996}
P.~Donnelly and T.G. Kurtz.
\newblock A countable representation of the {F}leming {V}iot measure-valued
  diffusion.
\newblock {\em Ann. Probab.} 24(2):698--742, 1996.

\bibitem{DonnellyKurtz1999}
P.~Donnelly and T.G. Kurtz.
\newblock Particle representations for measure-valued population models.
\newblock {\em Ann. Probab.} 27(1):166--205, 1999.

\bibitem{DrmotaEtAl2007}
M.~Drmota, A.~Iksanov, M.~M\"ohle, and U.~R\"osler.
\newblock {Asymptotic results concerning the total branch length of the
  Bolthausen$-$Sznitman coalescent}.
\newblock {\em Stochastic. Process. Appl.} 117(10):1404--1421, 2007.

\bibitem{Durrett2008}
R.~Durrett.
\newblock {\em Probability {M}odels for {DNA} {S}equence {E}volution}.
\newblock Springer, second edition, 2008.

\bibitem{Etheridge2001}
A.~Etheridge.
\newblock {\em An introduction to superprocesses}.
\newblock American Mathematical Society, 2001.

\bibitem{EtheridgePfaffelhuberWakolbinger2006}
A.~Etheridge, P.~Pfaffelhuber, and A.~Wakolbinger.
\newblock An approximate sampling formula under genetic hitchhiking.
\newblock {\em Ann. Appl. Probab.} 15:685--729, 2006.

\bibitem{EthierKurtz1986}
S.N. Ethier and T.~Kurtz.
\newblock {\em Markov {P}rocesses. {C}haracterization and {C}onvergence}.
\newblock John Wiley, New York, 1986.

\bibitem{Evans2000}
S.~Evans.
\newblock Kingman's coalescent as a random metric space.
\newblock In {\em Stochastic Models: Proceedings of the International
  Conference on Stochastic Models in Honor of Professor Donald A. Dawson,
  Ottawa, Canada, June 10-13, 1998 (L.G Gorostiza and B.G. Ivanoff eds.)},
  Canad. Math. Soc., 2000.

\bibitem{FuLi1993}
Y.-X. Fu and W.-H. Li.
\newblock Statistical tests of neutrality of mutations.
\newblock {\em Genetics} 133:693--709, 1993.

\bibitem{GrevenPfaffelhuberWinter2010}
A.~Greven, P.~Pfaffelhuber, and A.~Winter.
\newblock Tree-valued resampling dynamics. {M}artingale problems and
  applications.
\newblock {\em Submitted}, 2010.

\bibitem{Kingman1982a}
J.~F.~C. Kingman.
\newblock The coalescent.
\newblock {\em Stochastic Process. Appl.} 13(3):235--248, 1982.

\bibitem{Kingman1982}
J.~F.~C. Kingman.
\newblock On the genealogy of large populations.
\newblock {\em J. Appl. Probab.} 19A:27--43, 1982.

\bibitem{Moehle2006}
M.~M{\"o}hle.
\newblock On the number of segregating sites for populations with large family
  sizes.
\newblock {\em Adv. Appl. Probab.} 38:750--767, 2006.

\bibitem{PfaffelhuberWakolbinger2006}
P.~Pfaffelhuber and A.~Wakolbinger.
\newblock The process of most recent common ancestors in an evolving
  coalescent.
\newblock {\em Stochastic Process. Appl.} 116:1836--1859, 2006.

\bibitem{Pitman1999}
J.~Pitman.
\newblock Coalescents with multiple collisions.
\newblock {\em Ann. Prob.} 27(4):1870--1902, 1999.

\bibitem{RauchBarYam2004}
E.~M. Rauch and Y.~Bar-Yam.
\newblock Theory predicts the uneven distribution of genetic diversity within
  species.
\newblock {\em Nature} 431:449--452, 2004.

\bibitem{Saundersetal1984}
I.~W. Saunders, S.~Tavar\'e, and G.~A. Watterson.
\newblock On the genealogy of nested subsamples from a haploid population.
\newblock {\em Adv. Appl. Probab.} 16:471--491, 1984.

\bibitem{Tajima1990}
F.~Tajima.
\newblock Relationship between {DNA} polymorphism and fixation time.
\newblock {\em Genetics} 125:447--454, 1990.

\bibitem{Tavare1984}
S.~Tavar\'e.
\newblock Line-of-descent and genealogical processes and their applications in
  population genetics models.
\newblock {\em Theor. Pop. Biol.} 26:119--164, 1984.

\bibitem{Tavare2004}
S.~Tavar\'e.
\newblock {\em Ancestral Inference in Population Genetics, in: Lectures on
  Probability and Statistics 1Ð188, in: Lecture Notes in Mathematics, vol.
  1837}.
\newblock Springer, 2004.

\bibitem{Wakeley2008}
J.~Wakeley.
\newblock {\em Coalescent Theory: An Introduction}.
\newblock Roberts \& Company, 2008.

\bibitem{Watterson1982}
G.A. Watterson.
\newblock Mutant substitutions at linked nucleotide sites.
\newblock {\em Adv. Appl. Prob.} 14:166--205, 1982.

\bibitem{WiufHein1999}
C.~Wiuf and J.~Hein.
\newblock {{R}ecombination as a point process along sequences}.
\newblock {\em Theo. Pop. Biol.} 55:248--259, 1999.

\end{thebibliography}

\end{document}